\documentclass{amsart}
\usepackage{amsmath} 
\usepackage{amssymb}

\usepackage{mathrsfs}

\newtheorem{theorem}{Theorem}[section] 
\newtheorem{claim}[theorem]{Claim}

\theoremstyle{definition}
\newtheorem{definition}[theorem]{Definition}

\newtheorem{convention}[theorem]{Convention}

\newtheorem{question}[theorem]{Question}
\newtheorem{observation}[theorem]{Observation} 
\newtheorem{discussion}[theorem]{Discussion}
\newtheorem{conjecture}[theorem]{Conjecture}

\newtheorem{hypothesis}[theorem]{Hypothesis}

\theoremstyle{remark}
\newtheorem{remark}[theorem]{Remark}
\newtheorem{notation}[theorem]{Notation}

\newcommand{\rest}{{\restriction}}

\newcommand{\ord}{{\rm ord}} 
\newcommand{\rd}{{\rm rd}} 
\newcommand{\tr}{{\rm tr}} 
\newcommand{\md}{{\rm md}} 
\newcommand{\EC}{{\rm EC}} 
\newcommand{\PC}{{\rm PC}} 
\newcommand{\Min}{{\rm Min}} 
\newcommand{\Mod}{{\rm Mod}} 
\newcommand{\Dom}{{\rm Dom}} 
\newcommand{\Ord}{{\rm Ord}} 
\newcommand{\Levy}{{\rm Levy}} 
\newcommand{\Rang}{{\rm Rang}} 
\newcommand{\Th}{{\rm Th}} 
\newcommand{\iif}{{\rm if}} 
\newcommand{\tp}{{\rm tp}} 
\newcommand{\oor}{{\rm or}} 
\newcommand{\odd}{{\rm odd}} 
\newcommand{\fin}{{\rm fin}} 
\newcommand{\lex}{{\rm lex}}

\newcommand{\wilog}{{\rm without loss of generality}}

\newcommand{\then}{{\underline{then}}}
\newcommand{\when}{{\underline{when}}}
\newcommand{\Then}{{\underline{Then}}}
\newcommand{\If}{{\underline{if}}}
\newcommand{\Iff}{{\underline{iff}}}
\newcommand{\mn}{{\medskip\noindent}}
\newcommand{\sn}{{\smallskip\noindent}}
\newcommand{\bn}{{\bigskip\noindent}}

\newcommand{\cD}{{\mathscr D}}

\newcommand{\cF}{{\mathscr F}}

\newcommand{\bbL}{{\mathbb L}}
\newcommand{\cM}{{\mathscr M}}

\newcommand{\bbP}{{\mathbb P}}
\newcommand{\cP}{{\mathscr P}}

\newcommand{\bbQ}{{\mathbb Q}}

\newcommand{\cT}{{\mathscr T}}
 
\newcommand{\cU}{{\mathscr U}}

\newcommand{\cf}{{\rm cf}}

\newcount\skewfactor
\def\mathunderaccent#1#2 {\let\theaccent#1\skewfactor#2
\mathpalette\putaccentunder}
\def\putaccentunder#1#2{\oalign{$#1#2$\crcr\hidewidth
\vbox to.2ex{\hbox{$#1\skew\skewfactor\theaccent{}$}\vss}\hidewidth}}
\def\name{\mathunderaccent\tilde-3 }

\newenvironment{PROOF}[2][\proofname.]
   {\begin{proof}[#1]\renewcommand{\qedsymbol}{\textsquare$_{\rm #2}$}}
   {\end{proof}}

\begin{document}

\title {Dependent $T$ and Existence of limit models}
\author {Saharon Shelah}
\address{Einstein Institute of Mathematics\\
Edmond J. Safra Campus, Givat Ram\\
The Hebrew University of Jerusalem\\
Jerusalem, 91904, Israel\\
 and \\
 Department of Mathematics\\
 Hill Center - Busch Campus \\ 
 Rutgers, The State University of New Jersey \\
 110 Frelinghuysen Road \\
 Piscataway, NJ 08854-8019 USA}
\email{shelah@math.huji.ac.il}
\urladdr{http://shelah.logic.at}
\thanks{The author thanks Alice Leonhardt for the beautiful typing.
This research was supported by the Israel Science Foundation.  First
 version typed in 14/Nov/2005. Paper 877}



% First Typed - 05/May/25 \null\newline
%First done 29/Apr/05

%Previous version - 2015/May/8

\date{June 2, 2015}

\begin{abstract}
Does the class of linear orders have (one of the variants of) the so
called $(\lambda,\kappa)$-limit model?  It is necessarily unique, and
naturally assuming some instances of G.C.H. we get some positive,
i.e. existence results.  More generally, letting $T$ be a complete first order 
theory and for simplicity assume G.C.H., for 
regular $\lambda > \kappa > |T|$ does $T$ have
(variants of) a $(\lambda,\kappa)$-limit models, except for stable
$T$?  For some, yes, the theory of dense linear order, for some, no.
Moreover, for independent $T$ we get negative, i.e. non-existence 
results.  We deal more
with linear orders.
\end{abstract}

\maketitle
\numberwithin{equation}{section}
\setcounter{section}{-1}
\newpage

\section {Introduction} 

The first part of the introduction is intended for a general
mathematical reader.
Cantor proved that the structure ``the rationals as a linear order" 
is characterized up
to isomorphism by being ``a dense linear order with neither first nor
last element which is countable".  Hausdorff generalizes this as
follows.  
For transparency 
assume the G.C.H., the generalized
continuum hypothesis then for every cardinal $\lambda$ there is a
unique linear order $I$ of cardinality $\lambda^+$ which is
$\lambda^+$-dense (i.e. if $A < C$ are subsets of cardinality $\le
\lambda$ then for some $b \in I$ we have $A < b < C$) with neither
first nor last elements).  This canonical linear order is, in later
model theoretic notions, the unique saturated model of the theory
$T_{\text{\rm ord}} = \text{\rm Th}(\bbQ,<)$ of cardinality
$\lambda^+$ (also the universal homogeneous model); note $T_{\ord}$ is the
first order theory of the rational order.

Later Bjarni J\'onsson \cite{Jo56}, \cite{Jo60} introduced and proved
the existence of homogeneous-universal models in cardinality $\lambda$,
for a quite general class of structures.  Morley and Vaught \cite{MoVa62}
introduced the notion of saturated models 
and investigate such models (which are homogeneous
universal if we use elementary submodels instead of
substructures).  Saturated models become a central notion in model
theory.

The author in \cite{Sh:88} or \cite{Sh:88r} = \cite[Ch.I]{Sh:h},
introduce abstract elementary classes and there define some variants of
$(\lambda,\kappa)$-limit models which are again (like the
homogeneous universal ones) unique but for the
pair of cardinals $\lambda,\kappa$; note that 
for $\lambda = \kappa = \mu^+$ this
is the previous case.  So natural questions are: what about 
elementary classes, i.e. first order theories? and what about the 
class of linear orders?

By \cite{Sh:868} if $T$ is low enough (so called stable) there are
existence theorems but, e.g. the theory of linear order is not stable.

What are our main results?  First, a result meaningful also to one
with very little set theoretic background.  If $\lambda =
\lambda^{<\lambda}$, e.g. $\lambda = \mu^+ = 2^\mu$,
 then in addition to the unique (up to
isomorphisms) linear order which Hausdorff discovers, for $\kappa =
\aleph_0$ or just $\kappa \le \lambda$ which is a successor (or just a
so called regular) 
there is a $(\lambda,\kappa)$-limit linear order and it is unique up
to isomorphisms.  We can also have a characterization (as in the case
of Hausdorff), though not so elegant; see \S1 which do not require model
theoretic background, see Theorem \ref{1.1}.  There are stronger versions of
``$(\lambda,\kappa)$-limit models" (($\lambda,\kappa$)-superlimit) for
which we show non-existence, see \S3.

Second, in model theoretic terms this shows that having
$(\lambda,\kappa)$-limit model is satisfied by some (complete first
order) theories $T$ which are not stable; all this in \S1.  So does
every $T$ have such models?  In \S2 comes another major result of
this work: the answer in general, is no, e.g. for (the first order
theory) Peano arithmetic, see Theorem \ref{2.2}.  Moreover, there is a
reasonable natural sufficient condition: the theory $T$ is so called
dependent, this is Theorem \ref{2p.3}.

Those complementary results lead to the main conjectures arising from
this work on existence of $(\lambda,\kappa)$-limit models and to the
generic pair conjecture.  They essentially say that the above
mentioned sufficient condition, ``$T$ is dependent" is the right one,
each dealing with a variant of the question (the first: any relevant
$\kappa$, the second: the parallel for $\kappa = 2$).

The question can be rephrased (under G.C.H., restricting ourselves to
successor cardinality $\aleph_{\varepsilon +1}$) as follows: assume $\langle
M_\alpha:\alpha < \aleph_{\varepsilon + 2}\rangle$ is a
$\prec$-increasing continuous sequence of models of the first order
complete $T,\|M_\alpha\| = \aleph_{\varepsilon +1}$ and $M =
\bigcup\{M_\alpha:\alpha < \aleph_{\varepsilon +2}\}$ is saturated
(e.g. Hausdorff linear order of cardinality $\aleph_{\varepsilon +2},T =
T_{\ord}$).  Let $\bold n_{\aleph_{\varepsilon + 1}}(T) =
\Min\{|\{M_\alpha/\cong:\alpha \in E\}|:E$ a closed unbounded
subset of $\aleph_{\varepsilon +2}\}$.  Now
the existence of $(\aleph_{\varepsilon +1},\kappa)$-limit model
for every regular $\kappa < \aleph_{\varepsilon +2}$ 
implies $\bold n_T(\aleph_{\varepsilon +1}) =
|\varepsilon +2|$, in fact for some such $E$ for any $\kappa$ all the models
$\{M_\delta:\delta \in E$ has cofinality
$\kappa\}$ are pairwise isomorphic.  Our non-existence results give
$\bold n_{\aleph_{\varepsilon +1}}(T) = \aleph_{\varepsilon +2}$.

In the rest of the introduction we assume more background.
\bigskip

\noindent
\centerline {$* \qquad * \qquad *$}
\bigskip

\noindent
We continue \cite{Sh:868} and \cite{Sh:783}

The problem in \cite{Sh:868} is when does (a first order theory) $T$ have a 
model $M$ of cardinality $\lambda$ which is (one of the
variants of) a limit model for cofinality $\kappa$, in the cases not
covered by \cite[0.8]{Sh:868} (or \cite[3.3,3.2]{Sh:88}, 
\cite[3.6,3.5]{Sh:88r}).  
More accurately, there are some versions of limit models, ``$M$ is a
$(\lambda,\kappa)$-$x$-limit model of $T$" mainly
``$(\lambda,\kappa)$-i.md. limit", see Definition \ref{y.5}; (though
we deal with others, too) the most natural case to try is $\lambda =
\lambda^{< \lambda} > \kappa = \text{ cf}(\kappa) > |T|$.

Note that if $T$ has (any version of) a limit model of cardinality
$\lambda$ then there is a universal $M \in \Mod_\lambda(T)$.
Now we know that if $\lambda = 2^{< \lambda} > |T|$ then there is a
universal $M \in \Mod_\lambda(T)$ (see e.g. \cite{Ho93}).  But
for other cardinals it is ``hard to have a universal model",
see history \cite{KjSh:409} and \cite{Dj05}.  
E.g. if $T$ has the strict order property, then, by Kojman-Shelah
\cite{KjSh:409} there are ZFC non-existence results (a major case, for regular 
$\lambda$ is when 
$(\exists \mu)(\mu^+ < \lambda \wedge 2^\mu > \lambda)$.  In at least
one case, $\lambda = \aleph_1 < 2^{\aleph_0}$ consistently we do not
have a universal model, see \cite{Sh:100}.

Stable theories have limit models (in many cases); hence it is natural to ask:
\medskip

\noindent
\underline{Question 1}:  Assume
$\lambda = \lambda^{< \kappa} > \kappa > |T|$.  Does the existence of
a $(\lambda,\kappa)$-md.-limit model of $T$ imply $T$ is stable?

This is quite reasonable but in Theorem \ref{1.1} we 
find a counterexample, in fact, one everyone
knows about: the theory $T_{\ord}$ of
dense linear orders (see \ref{0.2}).  
This per se is a continuation of Hausdorff result,
revealing some canonical linear ordres.  Returning to the family of
elementary classes, i.e. first order theories, it is natural to ask:
\medskip

\noindent
\underline{Question 2}:  Does $T$ have a
$(\lambda,\kappa)$-i.md.-limit model whenever $\lambda = \lambda^{<
  \lambda} > \kappa + |T|$ for every unstable $T$?  

For non-existence results it is natural to look at $T$ dissimilar to
$T_{\ord}$.

\noindent
As $T_{\ord}$ is 
prototypical of dependent theories, it is natural to look for
independent theories.  A strong, explicit version of $T$ being
independent is having the strong independence property (see Definition
\ref{2.2.1}), e.g. Peano arithmetic has.  We prove that
for such $T$ there are no limit models (\ref{2.2}).  
But the strong independence property does not seem a good dividing
line.  The independence property is a good candidate for being a
meaningful dividing line.
\medskip

\noindent
\underline{Question 3}:  If $T$ is independent, does $T$ have a
$(\lambda,\kappa)$-i.md.-limit model (with $\lambda = 
\lambda^{< \lambda} > \kappa > |T|$)?

We work harder (than in \ref{2.2}) to prove (in \ref{2p.3}) the 
negative answer 
for every independent $T$ (for many cardinals), i.e.  with 
the independence property though a weaker version meaning we prove
non-existence of a stronger version of ``$(\lambda,\kappa)$-limit model".

This makes us
\begin{conjecture}
\label{0.x.7}  
Any dependent $T$ has $(\lambda,\kappa)$-i.md.-limit model.

Toward this end we intend to continue the investigation of types for dependent
$T$.

We shall also consider a property $\Pr_{\lambda,\kappa}(T)$ (and the
stronger $\Pr^2_{\lambda,\kappa}(T))$, see
Definition \ref{2.2.7}, which are relatives of ``there is no
$(\lambda,\kappa)$-$x$-limit model"; i.e. non-existence results for
independent $T$ holds for $\lambda = \lambda^{< \lambda} \ge \kappa =
\cf(\kappa),\lambda > |T|$.  For $\lambda > \kappa$ this
strengthens ``there is no $(\lambda,\kappa)$-i.md.-limit model".  But
$\lambda = \kappa$ is a new non-trivial case and it is also a
candidate to be ``an outside equivalent condition for $T$ being
dependent".
\end{conjecture}

\noindent
The most promising among the relatives (for having a dichotomy) 
is the following conjecture (the assumption $2^\lambda = \lambda^+$ is
just for simplicity).
\begin{conjecture}
\label{0.x.16}
\underline{The generic pair conjecture}  

Assume $\lambda = \lambda^{< \lambda} > |T|$ and
$2^\lambda = \lambda^+$ (for transparency) and $M_\alpha \in 
\EC_\lambda(T)$ is $\prec$-increasing continuous for $\alpha <
\lambda^+$ with $\bigcup\{M_\alpha:\alpha < \lambda^+\} \in 
\EC_{\lambda^+}(T)$ saturated. Then $T$ is dependent iff for some club
$E$ of $\lambda^+$ for all pairs $\alpha < \beta <\lambda^+$ from $E$ both
of cofinality $\lambda,(M_\beta,M_\alpha)$ has the same isomorphism
type (we denote this property of $T$ by $\Pr^2_\lambda(T)$),
see Definition \ref{2.2.7}).

Here we prove that for independent $T$, a strong version of the
conjecture holds.  

In \S2, we also prove the
parallel of what we say above.  In \S3 we prove that
$(\lambda,\kappa)$-superlimit models does not exist even for $T =
T_{\ord}$.   This work is continued in \cite{Sh:906},
\cite{Sh:900}, \cite{Sh:950}, \cite{Sh:F1124} and Kaplan-Lavi-Shelah 
\cite{KpLaSh:1055}.
\end{conjecture}
\bigskip

\noindent
\centerline {$* \qquad * \qquad *$}
\bn 
Now we define some versions of ``$M$ is a $(\lambda,S)$-$x$-limit model"
and for them ``$\bar M$ obeys a $(\lambda,T)$-$x$-function".

\begin{notation}
\label{y.1}  
1) Let $T$ denote a complete first order theory.

\noindent
2) Let $\tau_T = \tau(T),\tau_M = \tau(M)$ 
be the vocabulary of $T,M$ respectively.
\end{notation}

\begin{definition}
\label{y.2}  
1) For any $T$ let $\EC(T) = \{M:M$ a $\tau_T$-model of $T\}$.

\noindent
2) $\EC_\lambda(T) = \{M \in \EC(T):M$ is of cardinality $\lambda\}$ 
and $\EC_{\lambda,\kappa}(T) = \{M \in \EC_\lambda(T):M$ 
is $\kappa$-saturated$\}$.

\noindent
3) We say $M \in \EC(T)$ is $\lambda$-universal when every $N
\in \EC_\lambda(T)$ can be elementarily embedded into $M$.

\noindent
4) We say $M \in \EC(T)$ is universal when it is
$\lambda$-universal for $\lambda = \|M\|$.

\noindent
5) For $T \subseteq T'$ let

\[
\PC(T',T) = \{M \restriction \tau_T:M \text{ is model of } T'\}
\]

\[
\PC_\lambda(T',T) = \{M \in \PC(T',T):M \text{ is of cardinality } \lambda\}.
\]
\end{definition}

\begin{definition}
\label{y.3}  
Given $T$ and $M \in \EC_\lambda(T)$ we say that 
$M$ is a $(\lambda,\kappa)$-superlimit model when: $M$ is a 
$\lambda$-universal model of cardinality $\lambda$ and if 
$\delta < \lambda^+$ is a limit ordinal such that $\cf(\delta) = \kappa,
\langle M_\alpha:\alpha \le \delta \rangle$ is
$\prec$-increasing continuous, and $M_{\alpha +1}$ is isomorphic to $M$ for
every $\alpha < \delta$ \then \, $M_\delta$ is isomorphic to $M$.
\end{definition}

\begin{remark}
\label{y.4}  
We shall use:
\mn
\begin{enumerate}
\item[$(a)$]   $(\lambda,\kappa)$-i.md.-limit in \ref{1.1},
(existence for $T_{\ord}$)
\sn
\item[$(b)$]   $(\lambda,\kappa)$-wk-limit in \ref{2.2},
(non-existence from ``$T$ is strongly independent")
\sn
\item[$(c)$]  $(\lambda,\kappa)$-md.-limit in \ref{2p.3},
(non-existence for independent $T$)
\sn
\item[$(d)$]   $(\lambda,\kappa)$-i.st.-limit for $T_{\ord}$:
\ref{sl.21} and \ref{sl.2}(3), \ref{sl.2.3}(3), 
(on characterization) for $T_{\ord}$)
\sn
\item[$(e)$]   $(\lambda,\kappa)$-superlimit in \ref{nl.11} (non-existence).
\end{enumerate}
\end{remark}

Recall the definition of some versions of ``$(\lambda,\kappa)$-limit model".
\begin{convention}
\label{y.4d}
In this work let ``$M$ is $(\lambda,S)$-limit" mean ``$M$ is
$(\lambda,S)-\md$-limit, see Definition below; similarly for
$(\lambda,\kappa)$. 
\end{convention}

\begin{definition}
\label{y.5}  Let $\lambda$ be a cardinal $\ge |T|$.  
For parts 3) - 5) but not 6), for simplifying
the presentation we assume the axiom of global choice; alternatively
restrict yourself to models with universe an ordinal $\in
[\lambda,\lambda^+)$.  Below if $S = \{\delta <
\lambda^+:\cf(\delta) = \kappa\}$ then instead $(\lambda,S)$ we may
write $(\lambda,\kappa)$, this is the main case.

\noindent
1) Let $S \subseteq \lambda^+$ be stationary.  A model 
$M \in \EC_\lambda(T)$
is called $(\lambda,S)$-st-limit (or $S$-strongly limit or 
$(\lambda,S)$-strongly limit) \when \, for some function:  
$\bold F:\EC_\lambda(T)\rightarrow \EC_\lambda(T)$ we have:
\mn
\begin{enumerate}
\item[$(a)$]  for $N \in \EC_\lambda(T)$ 
we have $N \prec \bold F(N)$
\sn
\item[$(b)$]   if $\delta \in S$ is a limit ordinal
and $\langle M_i:i < \delta
\rangle$ is a $\prec$-increasing continuous sequence 
\footnote{No loss if we add $M_{i+1} \cong M$, so this simplifies the
demand on $\bold F$, i.e., only $\bold F(M)$ is required}
in $\EC_\lambda(T)$ obeying $\bold F$ which means
$i < \delta \Rightarrow \bold F(M_{i+1}) \prec M_{i+2}$, 
\then \, $M \cong \cup\{M_i:i < \delta\}$.
\end{enumerate}

\noindent
2) Let $S \subseteq \lambda^+$ be stationary.  $M \in 
\EC_\lambda(T)$ is called $(\lambda,S)$-nr-limit (or
$S$-normally limit, or may omit nr/normally) 
\when \, for some function $\bold F:\EC_\lambda(T) \rightarrow 
\EC_\lambda(T)$ we have:
\mn
\begin{enumerate}
\item[$(a)$]   for every $N \in \EC_\lambda(T)$ we have 
$N \prec \bold F(N)$ 
\sn
\item[$(b)$]    if $\langle M_i:i < \lambda^+ \rangle$ is
a $\prec$-increasing continuous sequence of members of 
$\EC_\lambda(T),\bold F(M_{i+1}) \prec M_{i+2}$ \then \,
for some closed unbounded \footnote{We can use a filter as a parameter}
subset $C$ of $\lambda^+$,

\[
[\delta \in S \cap C \Rightarrow M_\delta \cong M].
\]
\end{enumerate}
\mn
2A) $M \in \EC_\lambda(T)$ is $(\lambda,S)$-limit$^+$ \when \, if
$\langle M_\alpha:\alpha < \lambda^+\rangle$ is $\subseteq$-increasing and
continuous and $\alpha < \lambda^+ \Rightarrow 
M_{\alpha +1} \cong M$ \then \, for some club $E$ of
$\lambda$ we have $\alpha \in E \cap S \rightarrow M_\alpha \cong M$. Notice
that being a $(\lambda,S)$-limit$^+$ implies being a $(\lambda,S)$-nr-limit.

\noindent
3) We define ``$M$ is $(\lambda,S)$-wk-limit", ``$(\lambda,S)$-md-limit" 
like ``$(\lambda,S)$-nr-limit", ``$(\lambda,S)$-st-limit" 
respectively by demanding that the domain of
$\bold F$ is the family of $\prec$-increasing continuous
sequences of members of $\EC_\lambda(T)$ of length $< \lambda^+$ and
replacing ``$\bold F(M_{i+1}) \prec M_{i+2}$" by 
``$M_{i+1} \prec \bold F(\langle M_j:j \le i+1 \rangle) 
\prec M_{i+2}"$.  (They are also called $S$-weakly limit, $S$-medium
limit, respectively.) 

\noindent
3A) We replace ``limit" by ``limit$^-$" if
$``\bold F(M_{i+1}) \prec M_{i+2}",``M_{i+1} \prec
\bold F(\langle M_j:j \le i+1 \rangle) \prec M_{i+2}"$ are replaced 
by $``\bold F(M_i) \prec M_{i+1}",``M_i \prec 
\bold F(\langle M_j:j \le i \rangle) \prec M_{i+1}"$ respectively.

\noindent
4) If $S = \lambda^+$ then we may omit $S$ (in parts (3), (4), (5)).  

\noindent
5) For $\Theta \subseteq \{\mu:\mu \le \lambda$ and $\mu$ is 
regular$\},M$ is $(\lambda,\Theta)$-strongly limit \If \, $M$ is 
$\{\delta < \lambda^+:\cf(\delta) \in \Theta\}$-strongly limit
in the sense of part (1).  
Similarly for the other notions (where 
$\Theta \subseteq \{\mu:\mu \text{ regular } \le \lambda\}$ is
non-empty and $S_1 \subseteq \{\delta < \lambda^+:\cf(\delta) \in \Theta\}$
is a stationary subset of $\lambda^+$).  If we do not write $\lambda$
we mean $\lambda = \|M\|$.  

\noindent
6) We say that $M \in K_\lambda$ is $(\lambda,S)$-i.st-limit (or 
$S$-invariantly strong limit) \when \, in
part (3), $\bold F$ is just a subset of $\{(M,N)/\cong:M
\prec N$ are from $\EC_\lambda(T)\}$ and in clause (b) of part (3)
we replace ``$\bold F(M_{i+1}) \prec M_{i+2}"$ by ``$(\exists
N)(M_{i+1} \prec N \prec M_{i+2} \wedge ((M_{i+1},N)/\cong) \in 
\bold F)$".  But abusing notation we still 
write $N = \bold F(M)$ instead $((M,N)/ \cong) \in \bold F$.  
Similarly with the other notions, i.e., we use the isomorphism type of
$\bar M \char 94 \langle N\rangle$.  
\end{definition}

\begin{observation}
\label{y.5c}  
1) If $\bold F_1,\bold F_2$ are as above and
$\bold F_1(N) \prec \bold F_2(N)$ (or $\bold F_1(\bar M) \prec 
\bold F_2(\bar M))$
whenever defined then if $\bold F_1$ is a witness so is $\bold F_2$.

\noindent
2) All versions of limit models imply being a universal model in
$\EC_\lambda(T)$. 
\end{observation}
\bigskip

\noindent
3) \underline{Obvious implication diagram}:
\label{y.6}  
For stationary $S \subseteq S^{\lambda^+}_\kappa$ as in \ref{y.5}(7):

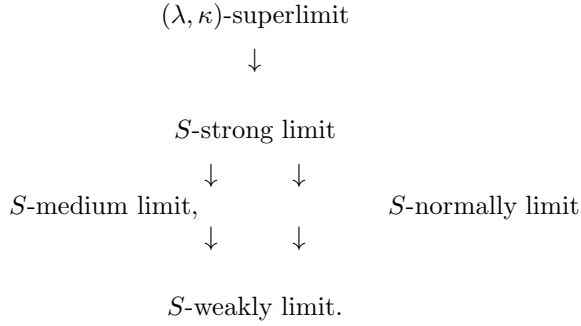
$$
(\lambda,\kappa)\text{-superlimit}
$$
\centerline {$\downarrow$}

$$
S \text{-strong limit}
$$

\centerline {$\downarrow$ \hskip25pt $\downarrow$}

\hskip75pt $S \text{-medium limit}, \qquad \qquad \qquad \quad 
S\text{-normally limit}$

\centerline {$\downarrow$ \hskip25pt $\downarrow$}
$$
S \text{-weakly limit}.
$$

\begin{claim}
\label{y10}
Assume $\lambda = \lambda^{< \kappa} \ge |T|$ and $\kappa$ is regular
and $M$ is a model of $T$ of cardinality $\lambda$.  \Then \, the
following conditions are equivalent (assuming the universal axiom of
choice \underline{or} restrict ourselves below to models with universe
$\subseteq \lambda^+$):
\mn
\begin{enumerate}
\item[$(a)$]  $M$ is $(\lambda,\kappa)$-md-limit
\sn
\item[$(b)$]  in the following game the isomorphism player has a
winning strategy.  A play last $\kappa$-moves, in the $i$-th move the
anti-isomorphism player chooses $M_\alpha \in \EC_\lambda(T)$ 
such that $\langle M_\beta:\beta \le \alpha\rangle$ is
$\prec$-increasing continuous and $\alpha = \beta + 1 \Rightarrow
M'_\beta \prec M_\alpha$ and the isomorphism player chooses
$M'_\alpha$ such that $M_\alpha \prec M'_\alpha \in 
\EC_\lambda(T)$.  The isomorphism player wins a play when 
$\bigcup\{M_\alpha:\alpha < \kappa\}$ is
isomorphic to $M$
\sn
\item[$(c)$]  there is a function $\bold F$ with domain $\{\bar M:\bar
  M$ a $\prec$-increasing continuous sequence of members of
  $\EC_\lambda(T)$ of length $< \kappa\}$ such that $i < \ell g(\bar
  M) \Rightarrow M_i \prec \bold F(\bar M)$ and if $\bar M = \langle
  M_i:i \le \kappa\rangle$ is $\prec$-increasing continuous sequence
  of members of $\EC_\lambda(T)$ and $i < \kappa \Rightarrow \bold
F(\bar M \rest (2i+2) \prec M_{2i+2}$ then $M_\kappa \cong M$
\sn
\item[$(d)$]  there is a function $\bold F$ such that: if $\langle
M_i:i \le \kappa\rangle$ is $\prec$-increasing continuous in
$\EC_\lambda(T)$ and for some sequence $\langle M'_i:i < \kappa\rangle$
we have $M_i \prec M'_i \in \EC_\lambda(T)$ and $i < \kappa
\Rightarrow M'_i \prec \bold F(\langle M_j:j \le i \rangle \char 94
\langle M'_i\rangle) = M_{i+1}$ (we say $\bar M$ obeys $\bold F$)
\then \, $\cup\{M_i:i < \kappa\} \cong M$
\sn
\item[$(e)$]  in $\bold V^{\Levy(\lambda^+,2^\lambda)}$ we
have:  if $\langle M_\alpha:\alpha < \lambda^+\rangle$ is 
$\prec$-increasing continuous, $M_\alpha \in \EC_\lambda(T)$,  
and $\bigcup\{M_\alpha:\alpha < \lambda^+\} \in
\EC_{\lambda,\lambda}(T)$ \then \, for some club $E$ of
$\lambda^+$ we have $\delta \in E \wedge \cf(\delta) =
\kappa \Rightarrow M_\delta \cong M$
\sn
\item[$(f)$]  like (e) for any $\lambda^+$-complete forcing notion
$\bbP$ such that $\Vdash_{\bbP} ``2^\lambda = \lambda^+"$.
\end{enumerate}
\end{claim}

\begin{PROOF}{\ref{y10}}
  As $(\EC_\lambda(T),\prec)$ has the JEP (joint
embedding property) and the amalgamation property this is
straightforward.

\noindent
E.G.
\medskip

\noindent
\underline{$(f) \Rightarrow (b)$}:

Let $\langle \name M_\alpha:\alpha < \lambda^+\rangle$ be a
$\bbP$-name of a $\prec$-increasing continuous sequence of members of
$\EC_\lambda(T)$ with union in $\EC_{\lambda,\lambda}(T)$ and $\name
E$ a $\bbP$-name of a club of $\lambda^+$ such that $\delta \in \name
E \wedge \cf(\delta) = \kappa \Rightarrow M_\delta \cong M$; clearly
it exists by clause (f) which we are assuming.  We now define a
strategy {\bf st} for the isomorphic player: together with choosing
$M'_\alpha$ the isomorphic player chooses
$(\gamma_\alpha,p_\alpha,h_\alpha)$ such that
\mn
\begin{enumerate}
\item[$\bullet_1$]  $M_\alpha \prec M'_\alpha \in \EC_\lambda(T)$ as
  demanded in (b)
\sn
\item[$\bullet_2$]  $p_\alpha \in \bbP$ and $\beta < \alpha
  \Rightarrow \bbP \models ``p_\beta \le p_\alpha"$
\sn
\item[$\bullet_3$]  $\gamma_\alpha < \lambda^+$ and $\beta < \alpha
  \Rightarrow \gamma_\beta < \gamma_\alpha$
\sn
\item[$\bullet_4$]  $p_\alpha \Vdash_{\bbP} ``h_\alpha$ is an
  isomorphism from $M'_\alpha$ onto $\name M_{\gamma_\alpha}"$
\sn
\item[$\bullet_5$]  if $\beta < \alpha$ then $h_\beta \subseteq
 h_\alpha$.
\end{enumerate}
\end{PROOF}

\noindent
Like \ref{y10} but for the invariant version we note
\begin{claim}
\label{y12}
For $M \in \text{\rm EC}_\lambda(T)$ the following are equivalent (and
seemingly stronger than the conditions in \ref{y10}): 
\mn
\begin{enumerate}
\item[$(a)'$]  $M$ is $(\lambda,\kappa)$-i.md-limit (that is invariantly
medium $(\lambda,\kappa)$-limit)
\sn
\item[$(d)'$]  there is a class $\bold F$ such that:
\sn
\begin{enumerate}
\item[$(\alpha)$]  $\bold F \subseteq \{\bar M:\bar M = \langle M_i:i \le
\alpha\rangle$ for some $\alpha \le \kappa$ is $\prec$-increasing
continuous, $\{M_i:i \le \alpha\} \subseteq \EC_\lambda(T)\}$ and
$\bold F$ is closed under isomorphisms
\sn
\item[$(\beta)$]  if $\bar M = \langle M_i:i \le \alpha\rangle \in
\bold F$ and $M_\alpha \prec M'_\alpha \in \EC_\lambda(T)$
then for some $M_{\alpha +1}$ we have $M'_\alpha \prec M_{\alpha +1}$
and $\bar M \char 94 \langle M_{\alpha +1}\rangle \in \bold F$
\sn
\item[$(\gamma)$]  for $\alpha$ limit $\langle M_i:i \le \alpha\rangle \in
\bold F$ iff $j < \alpha \Rightarrow \langle M_i:i \le j\rangle \in
\bold F$ and $M_\alpha = \cup\{M_i:i < \alpha\}$
\sn
\item[$(\delta)$]  if $\langle M_i:i \le \kappa\rangle \in \bold F$
then $M_\kappa \cong M$
\end{enumerate}
\sn
\item[$(d)''$] there is $\bold F$ such that:
\sn
\begin{enumerate}
\item[$(\alpha)$]  $\bullet \quad 
\bold F$ is a subset of $\{\bar M:\bar M = \langle M_i:i \le \alpha 
\rangle$ for some $\alpha < \kappa$ is 

\hskip30pt $\prec$-increasing continuous in $\EC_\lambda(T)\}$ 
\sn
\item[${{}}$]  $\bullet \quad$ [$\alpha$ odd $\Rightarrow M_\alpha 
\prec \bold F(\bar M) \in \EC_\lambda(T)$) so we can ignore the member

\hskip30pt $\alpha \rangle)$]
\sn
\item[$(\beta)$]  $\bullet \quad$ if $\bar M \in K$ has length $2
\alpha +1 < \kappa$ \then \, for some $M',\bar M \char 94 \langle (M')
\rangle \in \bold F$ 

\hskip30pt  and $M'$ is unique up to isomorphism, i.e. if
$\bar M^\ell =$ 

\hskip30pt  $\bar M \char 94 \langle M^\ell \rangle \in \cF$ for
$\ell=1,2$ then $\bar M^1,\bar M^2$ are isomorphic, 

\hskip30pt  so abusing notation we may write $M' = \bold F(\bar M)$
\sn
\item[${{}}$]  $\bullet \quad$ if $\bar M \in K$ has length $2
\alpha < \kappa$ and $M_{2 \alpha} \prec M' \in \EC_\lambda(T)$ then 

\hskip25pt $\bar M \char 94 \langle \bar M' \rangle \in \bold F$ and
we may write $\bold F(\bar M) = M_\alpha$
\sn
\item[$(\gamma),(\delta)$]   as in $(d)'$.
\end{enumerate}
\end{enumerate}
\end{claim}
\bigskip

\begin{definition}
\label{0.2}  
1) $T_{\ord}$ is the theory of dense linear 
order with neither first nor last element.

\noindent
2) $T_{\rd}$ is the theory of linear orders.
\end{definition}

\begin{definition}
\label{z15}
1) We say that $(C_1,C_2)$ is a cut of $M \in \EC(T_{\rd})$ when:
\mn
\begin{enumerate}
\item[$(a)$]  $C_1$ is an initial segment of $M$
\sn
\item[$(b)$]  $C_2$ is an end-segment of $M$
\sn
\item[$(c)$]  $C_1 \cap C_2 = \emptyset$
\sn
\item[$(d)$]  $C_1 \cup C_2 = M$.
\end{enumerate}
\mn
2) For a cut $(C_1,C_2)$ of $M$, let $\cf(C_1,C_2)$, the cofinality of
the cut $(C_1,C_2)$, be the pair $(\theta_1,\theta_2)$ when
\mn
\begin{enumerate}
\item[$(a)$]  $\theta_1$ is the cofinality of $C_1$, i.e. of $M
\restriction C_1$ (can be $0,1$ or a regular cardinal $\in
[\aleph_0,\lambda)$
\sn
\item[$(b)$]  $\theta_2$ is the cofinality of $C_2$ inverted.
\end{enumerate}
\end{definition}

\begin{definition}
\label{y17}
$\dot I(\lambda,T)$ is the number of $M \in \EC_\lambda(T)$ up to isomorphism.
\end{definition}

\begin{definition}
\label{y21}
1) Fixing $T,\varphi(\bar x,\bar y)$ is an independent formula when
for every $n$ and $M \models T$ for some $\bar a_\ell \in {}^{\ell
  g(\bar y)}M$ for $\ell < n$, for every $u \subseteq n,M \models
(\exists \bar x) \bigwedge\limits_{\ell < n} \varphi(\bar x,\bar
a_\ell)^{\iif(\ell \in u)}$.

\noindent
2) $T$ is independent iff some $\varphi(x,\bar y)$ is independent. 
\end{definition}

\begin{notation}
%\label{}
$\varphi^{\iif(\bold t)}$ is $\varphi$ if $\bold t$ is true or $1,\neg
  \varphi$ if $\bold t$ is false or 0.
\end{notation}

\begin{definition}
\label{y23}
1) $\lambda^{<\kappa>_{\tr}}$ is $\sup\{|\cT \cap {}^\kappa \lambda|:\cT
   \subseteq {}^{\kappa \ge}\lambda$ is closed under initial segments
   and $\varepsilon < \kappa \Rightarrow |\cT \cap {}^\varepsilon
   \lambda| \le \lambda\}$.

\noindent
2) $\bold U_\kappa(\lambda) = \min\{|\cF|:\cF$ is a set of functions from
   $\kappa$ to $\lambda$ such that $f \ne g \in \cF \Rightarrow \kappa >
   \{i < \kappa:f(i) = g(i)\}|$.
\end{definition}
\newpage

\section {Dense linear order has medium limit models} 

\begin{theorem}
\label{1.1}  
Assume $\lambda = \lambda^{< \lambda}
>\kappa = \cf(\kappa)$.  \Then \, $T_{\ord}$ 
 has an invariantly medium $(\lambda,\kappa)$-limit model.
\end{theorem}

\begin{remark}
\label{a5}
1) We use condition $(d)'$ from \ref{y12}, we may use it
as a definition.

\noindent
2) So a model of $T_{\text{ord}}$ is a dense linear order with neither
first nor last element and $\prec$ for models of $T_{\ord}$ is
just $\subseteq$ and saturated means $\lambda$-dense for models of
$T_{\ord}$ of cardinality $\lambda$.

\noindent
3) We actually prove a result with $\bold F$ of a simple kind, dealing
with $\bold F$ acting on pairs of models, $\cup\{M_i:i < \kappa\}$ is
isomorphic to the $(\lambda,\kappa)$-i.md.-limit model when $\langle
M_i:i < \kappa\rangle$ is $\prec$-increasing continuous sequence of
linear orders such that for any $i_1 < i_2 < \kappa$ for some $i_3 \in
(i_2,\kappa)$ we have $\bold F(M_{i_1},M_{i_3}) \prec M_{i_3+1}$.

\noindent
4) On cuts and their cofinalities \ref{z15}.
\end{remark}

\begin{remark}  
\label{a8}
Concerning $\circledast^\kappa_{\bar\mu}$ in the beginning of the
proof of \ref{1.1}.

\noindent
1) It is a characterization of the invariantly medium
   $(\lambda,\kappa)$-model.  We shall return to this in \S3.

\noindent
2) Concerning the clause inside $\circledast^\kappa_\mu$ in the proof
of \ref{1.1} note the following: Clause (f) almost implies clause (d).

Clause (f) implies (h)$_1$; why?  
use $A = \emptyset,B = M_i$.  Also clause (f) implies (h)$_2$; why? 
use $A = M_i,B = \emptyset$.

Lastly, clause (f) implies (i)$_1$; why?  use 
$A' = A,B' = \{c \in M_i:A < c\}$ and clause (f) implies (i)$_2$
similarly.

\noindent
2) Note that clause (f) is equivalent to (i)$_1$ + (i)$_2$.
\end{remark}

\begin{PROOF}{\ref{1.1}}
First we say that $\bar M$ is a fast
$(\lambda,\kappa)$-sequence (of models of $T_{\ord}$) when:
\mn
\begin{enumerate}
\item[$\circledast^\kappa_{\bar M}$]  $(a) \quad \bar M = \langle
 M_i:i \le \kappa \rangle$
\sn
\item[${{}}$]  $(b) \quad M_i$ is $\prec$-increasing continuous
\sn
\item[${{}}$]  $(c) \quad M_i$ is a model of $T_{\ord}$ of cardinality
$\lambda$
\sn
\item[${{}}$]  $(d) \quad M_i$ is saturated if $i$ is a non-limit ordinal
\sn
\item[${{}}$]  $(e) \quad$ if $i < \kappa$ and 
$a \in M_{i+1} \backslash M_i$ then
$M_{i+1} \restriction \{b \in M_{i+1} \backslash M_i:(\forall c \in M_i)$

\hskip25pt $[(c < b) \equiv (c < a)]\}$ is a saturated model of
$T_{\text{\rm ord}}$ of cardinality $\lambda$
\sn
\item[${{}}$]  $(f) \quad$ if $i < \kappa$ and $A,B \subseteq M_i$
and $A < B$ (i.e. $(\forall a \in A)(\forall b \in B)(a <_{M_i} b))$

\hskip25pt and $A$ or $B$ has 
cardinality $< \lambda$, \then \, for some $c \in M_{i+1} \backslash M_i$ 

\hskip25pt we have $A < c < B$; this includes $A,B$ singletons but it
is enough 

\hskip25pt to have this when $c \in M_i \Rightarrow \neg(A < c < B)$;
note that we say

\hskip25pt ``$A$ \underline{or} $B \ldots$"
\sn
\item[${{}}$]  $(g)_1 \quad$ if $i < j < \kappa$ and $a \in M_j
\backslash M_i$, \then \, for some $d \in M_{j+1} \backslash M_j$ we
have
\sn
\begin{enumerate}
\item[${{}}$]  $\quad (\alpha) \quad$ if $b \in M_i$ and $b <_{M_j} a$
then $b <_{M_{j+1}} d$
\sn
\item[${{}}$]  $\quad (\beta) \quad$ if $c \in M_j$ and $(\forall b
\in M_i)(b <_{M_j} a \Rightarrow b <_{M_j} c)$ then $d <_{M_{j+1}} c$
\end{enumerate}
\sn
\item[${{}}$]   $(g)_2 \quad$ if $i < j < \kappa$ and $a \in M_j
\backslash M_i$ \then \, for some $d \in M_{j+1} \backslash M_j$ we have
\sn
\begin{enumerate}
\item[${{}}$]  $\quad (\alpha) \quad$ if $b \in M_i$ and $a <_{M_j} b$
then $d <_{M_{j+1}} b$
\sn
\item[${{}}$]  $\quad (\beta) \quad$ if $c \in M_j$ and $(\forall b
\in M_i)(a <_{M_j} b \Rightarrow c <_{M_j} b)$ then $c <_{M_{j+1}} d$
\end{enumerate}
\sn
\item[${{}}$]  $(h)_1 \quad$ for $i < \kappa$ there is $b \in
M_{i+1} \backslash M_i$ such that $a \in M_i \Rightarrow a <_{M_{i+1}} b$
\sn
\item[${{}}$]   $(h)_2 \quad$ for $i < \kappa$ there is $b \in
M_{i+1} \backslash M_i$ such that $a \in M_i \Rightarrow b <_{M_{i+1}} a$
\sn
\item[${{}}$]   $(i)_1 \quad$ if $A \subseteq M_i,i < \kappa$ and
$|A| < \lambda$ \then \, for some $c \in M_{i+1} \backslash M_i$ we
have 

\hskip25pt $(\forall d \in M_i)(d <_{M_{i+1}} c \leftrightarrow (\exists a \in
A)(d \le_{M_i} a))$
\sn
\item[${{}}$]  $(i)_2 \quad$ if $A \subseteq M_i,i < \kappa$ and
$|A| < \lambda$ \then \, for some $c \in M_{i+1} \backslash M_i$ we
have 

\hskip25pt $(\forall d \in M_i)(c <_{M_{i+1}} d \leftrightarrow (\forall a \in
A)(a \le_{M_i} d))$
\sn
\item[${{}}$]   $(j) \quad$ if $i < \kappa$ and $a <_{M_i} b$
\then \, for some $c \in (a,b)_{M_{i+1}} \backslash M_i$ 

\hskip25pt the orders $M_i \restriction \{d \in M_i:d <_{M_{i+1}} c\}$
and the inverse of 

\hskip25pt $M_i \restriction \{d \in M_i:c <_{M_{i+1}} d\}$
have cofinality $\lambda$.
\end{enumerate}
\mn
Clearly:
\mn
\begin{enumerate}
\item[$\boxtimes$]  it is enough to prove $\boxtimes_1 + \boxtimes_2$ where
\sn
\begin{enumerate}
\item[$\boxtimes_1$]   there is $\bold F$ such that
\sn
\item[${{}}$]  $(\alpha) \quad \Dom(\bold F) = \{\bar M$: for some
  $\alpha \le \kappa,\bar M=\langle M_i:i \le \alpha \rangle$ is 

\hskip33pt $\prec$-increasing continuous, 
$M_i \in \EC_\lambda(T_{\ord})\}$
\sn
\item[${{}}$]  $(\beta) \quad$ for $\bar M = \langle M_i:i \le \alpha
  \rangle \in \bold F$
\sn
\item[${{}}$]  $\qquad \bullet \quad$ if $\alpha$ is odd then
  $M_\alpha \prec \bold F(\bar M) \in \EC_\lambda(T_{\Ord})$ and

\hskip33pt  $\bar M \char 94 \langle \bold F(\bar M)\rangle \in \bold F$
\sn
\item[${{}}$]  $\qquad \bullet \quad$ if $\alpha$ is even and
$M_\alpha \prec M' \in \EC_\lambda(T_{\Ord})$ then
 $M_\alpha = \bold F(\bar M)$

\hskip33pt and $\bar M \char 94 \langle M'\rangle \in \bold F$
\sn
\item[${{}}$]   $(\gamma) \quad \bold F$ is invariant, i.e. if
$\bar M_1 \cong \bar M_2$ then $(\bar M_1,\bold F(M_1)) \cong
(M_2,\bold F(\bar M_2))$
\sn
\item[${{}}$]  $(\delta) \quad$ if $\bar M = \langle
M_i:i \le \kappa \rangle$ is an $\prec$-increasing continuous sequence

\hskip33pt of members of $\EC_\lambda(T_{\ord})$ belonging to, i.e.
$\bold F$ \then \, $\circledast^\kappa_{\bar M}$
\end{enumerate}
\sn
\item[$\boxtimes_2$]  if $\circledast^\kappa_{\bar M^1}$ and
$\circledast^\kappa_{\bar M^2}$ then $M^1_\kappa \cong M^2_\kappa$.
\end{enumerate}
\bigskip

\noindent
\underline{Why is clause $\boxtimes_1$ true?}:

How do we choose ${\bold F}$?

Reading the definition of $\circledast^\kappa_{\bar M}$ this should be
clear: all our demands on $M_{j+1}$ when we are given 
$\langle M_i:i \le j \rangle$ and $M'_j$ can be fulfilled.  
We first choose ${\cP}_{\langle M_i:i \le j \rangle} = 
\{(A,B):(A,B)$ a cut of $M_j$ such that 
$A$ has cofinality $< \lambda$
\underline{or} the inverse of $B$ has cofinality $< \lambda$
\underline{or} for some $i < j$ and $a \in M_j \backslash M_i$ the set
$\{b \in M_i:b <_{M_j} a\}$ is unbounded in $A$ \underline{or} for
some $i < j$ and $a \in M_j \backslash M_i$ the set $\{b \in M_i:a
<_{M_j} b\}$ is unbounded from below in the set $B\}$.

Second, choose $M_{j+1} = \bold F(\langle M_i:i \le j\rangle \char 94
\langle M'_j\rangle)$ such that $M'_j \prec M_{j+1}$ and any cut $(A,B)
\in \cP_{\langle M_i:i \le j \rangle}$ is realized in 
$M_{j+1}$ and for each $c \in
M_{j+1} \backslash M_j$ we have $M_{j+1} \rest \{a \in M_{j+2}
\backslash M_j:a,c$ realize the same cut of $M_j\}$ is a saturated
model of $T_{\ord}$.

Having chosen $\bold F$, clauses $(\alpha),(\beta),(\gamma)$ of
$\boxtimes_1$ follow and clause $(\delta)$ follows too.
\bigskip

\noindent
\underline{Why is clause $\boxtimes_2$ true?}:

Suppose $\circledast^\kappa_{\langle M^\ell_i:i \le \kappa \rangle}$ for
$\ell=1,2$.

For $\ell=1,2$ let 

\begin{equation*}
\begin{array}{clcr}
Y_\ell = \{a \in M^\ell_\kappa \backslash M^\ell_0:&\text{ for every } A
\subseteq M^\ell_0 \text{ of cardinality } < \lambda \\
  &\text{ we have } A < a \Rightarrow (\exists b \in M^\ell_0)
(A < b < a) \text{ and} \\
  &a < A \Rightarrow (\exists b \in M^\ell_0)(a < b < A)\}
\end{array}
\end{equation*}

\[
E_\ell = \{(a,b):a,b \in M^\ell_\kappa \backslash M^\ell_0 
\text{ and } (\forall c \in M_0)(c < a \equiv c < b)\}.
\]
\mn
Now $E_\ell$ is an equivalence relation on $M^\ell_\kappa
\backslash M^\ell_0$ and $Y_\ell$ is a union of some equivalence
classes of $E_\ell$.  Let $Z_\ell \subseteq Y_\ell$ be a set of
representatives of $E_\ell \restriction Y_\ell$.  Now we define
$N_\ell$: it is the model with universe $M^\ell_0 \cup Z_\ell$, the
relation $<^{N_\ell} = <^{M^\ell_\kappa} \restriction (M^\ell_0 \cup
Z_\ell)$ and the relation $P^{N_\ell} = \{a:a \in M^\ell_0\}$.

Now it is easy to check that $N_\ell$ has first and last
elements both from $N_\ell \backslash P^{N_\ell}$ and is dense.
Also if $A,B \subseteq N_\ell$ have
cardinality $< \lambda$ and $A < B$ then we can find $a',a''$ such
that $A <_{N_\ell} \{a',a''\} <_{N_\ell} B$ and $a' \in P^{N_\ell},a''
\in N_\ell \backslash P^{N_\ell}$.  Hence $N_\ell$ is a saturated
model (not of $T_{\ord}$ but of a variant).   
So easily $N_1,N_2$ are isomorphic and let $g_0$ be 
such an isomorphism and $f_0 = g_0 \restriction M^1_0$.

Now
\mn
\begin{enumerate}
\item[$(*)_0$]   $f_0$ induces a mapping $\hat f_0$ from the class of
$E_1$-equivalence classes onto the class of $E_2$-equivalence classes.
\end{enumerate}
\mn
[Why?  Check the cases.]

Now we have to separately deal with each case of $M^1_\kappa \restriction
(a_1/E_1),M^2_\kappa \restriction (a_2/E_2)$ where $\hat f_0(a_1/E_1)
= a_2/E_2$.  But this is similar to the original problem, i.e., choose
$i < \kappa$ large enough such that $(a_1/E_1) \cap M^1_i \ne
\emptyset$ and $(a_2/E_2) \cap M^2_i \ne \emptyset$.  It is not hard
to understand that we can continue and in the end we exhaust the
models, but we shall elaborate; 
\wilog \, $M^1_\kappa \cap M^2_\kappa =\emptyset$.
For a set $A \subseteq M^\ell_\kappa$ we define
\mn
\begin{enumerate}
\item[$(*)_1$]   $E^\ell_A := \{(a,b):a,b \in M^\ell_\kappa
\backslash A$ and $(\forall c \in A)(a <_{M^\ell_\kappa} c \equiv b
<_{M^\ell_\kappa} c)\}$.
\end{enumerate}
\mn
Note
\mn
\begin{enumerate}
\item[$(*)_2$]   $E^\ell_A$ is an equivalence relation on
$M^\ell_\kappa \backslash A$.
\end{enumerate}
\mn
Define
\mn
\begin{enumerate}
\item[$(*)_3$]   $Y^\ell_A$ is the set $\{a \in M^\ell_\kappa
\backslash A$: the cut that $a$ induces on $A$ has cofinality
$(\lambda,\lambda)\}$.
\end{enumerate}
\mn
So
\mn
\begin{enumerate}
\item[$(*)_4$]   $Y^\ell_A$ is a subset $M^\ell_\kappa \backslash
A$ closed under $E^\ell_A$.
\end{enumerate}
\mn
Define
\mn
\begin{enumerate}
\item[$(*)_5$]  We say that $A \subseteq M^\ell_\kappa$ is $\ell$-nice
  \when \, for every $a \in M^\ell_\kappa \backslash A$,
for some $i=i_\ell(a,A) = i_\ell(a/E^\ell_A) < \kappa$ we have
\sn
\begin{enumerate}
\item[${{}}$]  $(\alpha) \quad$ the set $a/E^\ell_A$ is disjoint
to $M^\ell_i$ but not to $M^\ell_{i+1}$
\sn
\item[${{}}$]  $(\beta) \quad$ the set $\{b \in A:b <_{M^\ell_\kappa} a$
and $b \in M^\ell_i\}$ is unbounded in 

\hskip35pt $\{b \in A:b <_{M^\ell_\kappa} a\}$
\sn
\item[${{}}$]   $(\gamma) \quad$ the set $\{b \in A:a <_{M^\ell_\kappa} b$
and $b \in M^\ell_i\}$ is unbounded from 

\hskip35pt below in $\{b \in A:a <_{M^\ell_\kappa} b\}$
\end{enumerate}
\sn
\item[$(*)_6$]  in $(*)_5$, $i_\ell(a,A)$ is uniquely defined by
$(a,A)$, actually just by $a/E^\ell_A$
\sn
\item[$(*)_7$]   if $\delta < \kappa$ is a limit ordinal, $\ell \in
\{1,2\}$ and $\langle A_\alpha:\alpha < \delta\rangle$ is an
$\subseteq$-increasing sequence of $\ell$-nice sets such that $[\alpha
< \beta < \delta \wedge a \in M^\ell_\kappa \backslash A^\ell_\beta
\Rightarrow i(a,A_\alpha) < i(a,A_\beta)]$ then $A_\delta =:
\cup\{A_\alpha:\alpha < \delta\}$ is an $\ell$-nice set.
\end{enumerate}
\mn
[Why?  Trivially $A_\delta \subseteq M^\ell_\kappa$, so let $a \in
M^\ell_\kappa \backslash A_\delta$ then for each $\alpha < \delta$ we
have $a \in M^\ell_\kappa \backslash A_\alpha$ hence $i_\ell(a,A_\alpha) <
\kappa$ is well defined and it is $\le$-increasing with $\alpha$
because $(a/E^\ell_{A_\beta}) \subseteq (a/E^\ell_\alpha)$ by
clause $(\alpha)$ of $(*)_5$.  

Recall that $\langle i_\ell(a,A_\alpha):\alpha < \delta\rangle$ is
not eventually constant. 
We claim $i(*) = \bigcup\{i_\ell(a,A_\alpha):\alpha < \delta\}$ is as
required.
First of all, as the union of an $\le$-increasing not eventually
constant sequence of length $\delta
< \kappa$ of ordinals $< \kappa$ it is an ordinal $< \kappa$, in fact
a limit ordinal $< \kappa$.

Clearly, $a/E^\ell_{A_\delta}$ is the intersection of the 
$\subseteq$-decreasing
sequence $\langle a/E^\ell_{A_\alpha}:\alpha < \delta\rangle$.  Now if
$i<i(*)$ then for some $\alpha < \delta$ we have $i \le
i_\ell(a,A_\alpha)$ hence $a/E^\ell_{A_\alpha}$ is disjoint to $M^\ell_i$
hence $a/E^\ell_{A_\delta} \subseteq a/E^\ell_{A_\alpha}$ is disjoint
to $M_i$.  As this holds for every $i<i(*)$ it follows that also
$\bigcup\{M^\ell_i:i < i(*)\}$ is disjoint to $a/E^\ell_{A_\delta}$, but
$\bigcup\{M^\ell_i:i<i(*)\} = M^\ell_{i(*)}$ 
because $i(*)$ is a limit ordinal.  So
really $(a/E^\ell_{A_\delta}) \cap M^\ell_{i(*)} = \emptyset$.

It is also clear that $(\{b \in M^\ell_{i(*)}:b <_{M^\ell_\kappa}
a\},\{b \in M^\ell_{i(*)}:a <_{M^\ell_\kappa} b\})$ is a cut of
$M^\ell_{i(*)}$ whose cofinality $(\lambda_1,\lambda_2)$ is not equal
to $(\lambda,\lambda)$, hence by clauses $(i)_1 + (i_2)$ of
$\circledast^\kappa_{\bar M^\ell}$ we have $(a/E^\ell_{A_i}) \cap
M_{i(*)+1} \ne \emptyset$ so $i(*)$ satisfies demand $(\alpha)$ from
$(*)_5$ on $i(a,A_\delta)$.  The other two clauses should be clear, too.]

Define
\mn
\begin{enumerate}
\item[$(*)_8$]   ${\cF}$ is the set of $f$ such that
\begin{enumerate}
\item[${{}}$]   $(a) \quad$ for some 1-nice $A_1 \subseteq
M^1_\kappa$ and 2-nice set $A_2 \subseteq M^2_\kappa,f$ is an

\hskip25pt  isomorphism from the linear order 
$M^1_\kappa \restriction A_1$ onto 

\hskip25pt the linear order $M^2_\kappa \restriction A_2$
\sn
\item[${{}}$]  $(b) \quad$ for every $a_1 \in M^1_\kappa
\backslash A_1$ there is $a_2 \in M^2_\kappa \backslash A_2$ such that
$f$ maps 

\hskip35pt $\{b \in A_1:b < a_1\}$ onto $\{b \in A_2:b < a_2\}$; it
follows that

\hskip35pt  $a_1 \in Y^1_A$ iff $a_2 \in Y^2_A$
\sn
\item[${{}}$]  $(c) \quad$  for every $a_2 \in M^2_\kappa
\backslash A_2$ for some $a_1 \in M^1_\kappa \backslash A_1$ the
conclusion 

\hskip35pt of clause (b) holds.
\end{enumerate}
\end{enumerate}

Define
\mn
\begin{enumerate}
\item[$(*)_9$]   $<_*$ is the following two-place relation of ${\cF}:
f <_* f'$ \Iff \, $(f,f' \in {\cF}$ and)
\sn
\begin{enumerate}
\item[${{}}$]  $(a) \quad f \subseteq f'$
\sn
\item[${{}}$]  $(b) \quad$ if $a_1 \in M^1_\kappa
\backslash \Dom(f')$ then $i_1(a_1/E^1_{\Dom(f')}) > i_1(a_1/E^1_{\Dom(f)})$
\sn
\item[${{}}$]  $(c) \quad$ if $a_2 \in M^2_\kappa \backslash
\Rang(f')$ then $i_2(a_2/E^2_{\Rang(f')}) >
i_2(a_2/E^2_{\Rang(f)}$
\sn
\item[${{}}$]   $(d) \quad$ if $a \in M^1_\kappa \backslash
\Dom(f')$ then there are $b,c \in (a/E^1_{\Dom(f)}) \cap \Dom(f')$ 

\hskip25pt such that $b <_{M^1_\kappa} a <_{M^1_\kappa} c$
\sn
\item[${{}}$]   $(e) \quad$ if $a \in M^2_\kappa
\backslash \Rang(f')$ then there are $b,c \in (a/E^2_{\Rang(f)}) \cap$

\hskip25pt $\Rang(f')$ such that $b <_{M^2_\kappa} a <_{M^2_\kappa} c$.
\end{enumerate}
\end{enumerate}

Note
\mn
\begin{enumerate}
\item[$(*)_{10}$]  $({\cF},<_*)$ is a non-empty partial order.
\end{enumerate}
\mn
[Why?  We have in $(*)_0$ above proved that there is an isomorphism from
$M^1_0$ onto $M^2_0$ which belongs to ${\cF}$.  Being a partial order
is obvious.]
\mn
\begin{enumerate}
\item[$(*)_{11}$]   if $\delta < \kappa$ is a limit ordinal and
$\langle f_\alpha:\alpha < \delta\rangle$ is a $<_*$-increasing
sequence in ${\cF}$, \then \, $f_\delta := \bigcup\{f_\alpha:\alpha
< \delta\}$ belongs to ${\cF}$ and $\alpha < \delta \Rightarrow
f_\alpha <_* f_\delta$.
\end{enumerate}
\mn
[Why?  Clearly $f_\delta$ is an isomorphism from the linear order
$M^1_\kappa \restriction A^1_\delta$ where $A^1_\delta =:
\bigcup\{\Dom(f_\alpha):\alpha <\delta\}$ onto the linear order $M^2_\kappa
\restriction A^2_\delta$ where $A^2_\delta =:
\bigcup\{\Rang(f_\alpha):\alpha < \delta\}$.  Now $\Dom(f_\delta) =
\bigcup\{\Dom(f_\alpha):\alpha < \delta\}$ is 1-nice by $(*)_7$
recalling clause (a) of $(*)_9$ and
similarly $\Rang(f_\delta) = \bigcup\{\Rang(f_\alpha):\alpha <
\delta\}$ is $2$-nice.  So from the demands for ``$f_\delta \in
{\cF}$" in $(*)_8$, clause (a) holds. Concerning clause (b) there,
let $a_1 \in M^1_\kappa \backslash \Dom(f_\delta)$.  For each
$\alpha < \delta$ by $(*)_{10}(d)$ applied to $f_\alpha <_* f_{\alpha
+1}$ there is a pair $(b_\alpha,c_\alpha)$ satisfying 
$b_\alpha,c_\alpha \in (a/E^1_{\Dom(f_\alpha)}) \cap 
\Dom(f_{\alpha +1})$ such that $b_\alpha <_{M^1_\kappa} a_1
<_{M^1_\kappa} c_\alpha$.  Note that as $b_\alpha,c_\alpha \in
(a/E^1_{\Dom(f_\alpha)})$ necessarily $b_\alpha,c_\alpha
\notin \Dom(f_\alpha)$ and clearly $d \in \Dom(f_\alpha)
\Rightarrow (d < b_\alpha \equiv d < c_\alpha)$.  Hence $\langle
b_\alpha:\alpha < \delta\rangle$ is increasing, $\langle
c_\alpha:\alpha < \delta\rangle$ is decreasing, and: $d \in 
\Dom(f_\delta)$ implies that for some $\alpha < \delta,d \in
\Dom(f_\alpha)$ hence for every $\beta < \delta$
large enough $d < b_\beta \equiv d < c_\beta$.  Recall
$b_\alpha,c_\alpha \in \Dom(f_{\alpha +1}) \backslash 
\Dom(f_\alpha)$ so $\langle i_1(b_\alpha,\Dom(f_\beta)):
\beta \le\alpha\rangle$ is
increasing, $i_1(b_\alpha,\Dom(f_\beta)) = i_1(c_\alpha,
\Dom(f_\beta))$.  So $(\{b_\alpha;\alpha < \delta\},\{c_\alpha:\alpha <
\delta\})$ determine the cut $a_1$ induces on $\Dom(f_\delta)$
and they are $\subseteq M^1_{i_1(a,\Dom(f_\delta))}$.  
Now $(\{f_{\alpha +1}(b_\alpha):\alpha < \delta\},
(\{f_{\alpha +1}(c_\alpha):\alpha <
\delta\})$, behave similarly in $M^2_\kappa$ and they induce a cut of
$M_i,i = \bigcup\{i_2(f_{\alpha +1}(b_\alpha),\Rang(f_\alpha)):\alpha
< \delta\}$ which is realized by some $a_2 \in M_{i+1}$ by clause (f) of
$\circledast^2_{\bar M^2}$.  Now $a_2$ is as required. 

Clause (c) is proved similarly using $(*)_{10}(e)$.]
\mn
\begin{enumerate}
\item[$(*)_{12}$]   if $\langle f_\alpha:\alpha < \kappa\rangle$ is an
$<_*$-increasing sequence in ${\cF}$ then $f_\kappa :=
\bigcup\{f_\alpha:\alpha < \kappa\}$ is an isomorphism from $M^1_\kappa$
onto $M^2_\kappa$.
\end{enumerate}
\mn
[Why?  Toward contradiction first assume
$\Dom(f_\kappa) \subset M^1_\kappa$ so choose $a_1 \in
M^1_\kappa \backslash \Dom(f_\kappa)$ hence $\langle
i_1(a_1/E^1_{\Dom(f_\alpha)}):\alpha < \kappa \rangle$ is a
(strictly) increasing sequence of ordinals $< \kappa$, hence its sup is
$\kappa$.  Now for $\alpha < \kappa,a_1 \in
(a_1/D^1_{\Dom(f_\alpha)})$ but $a_1/E^1_{\Dom(f_\alpha)}$ is disjoint to 
$M^1_{i_1(a_2,\Dom(f_\alpha))}$.  Hence
$a_1 \notin \cup\{M^1_{i_1(a_2,\Dom(f_\alpha))}:\alpha < \kappa\} 
= \{M_\beta:\beta < \kappa\} = M_\kappa$
which is impossible.  Similarly $\Rang(f_\kappa) \subset
M^2_\kappa$ leads to contradiction, so we are done.]
\mn
\begin{enumerate}
\item[$(*)_{13}$]   for every $f \in {\cF}$ there is $f'$
such that $f <_* f' \in {\cF}$.
\end{enumerate}
\mn
[Why?  Let $\langle a^1_t:t \in I\rangle$ be a set of representatives
of $(M^1_\kappa \backslash \Dom(f))/E_{\Dom(f)}$.  For $t \in I$ choose
$a^2_t \in M^2_\kappa \backslash \Rang(f)$ such that $f$ maps
$\{b \in \Dom(f):b < a^1_t\}$ onto $\{b \in \Rang(f):b <
a^2_t\}$ and let $i_{1,t} := i_1(a^1_t/E^1_{\Dom(f)}),i_{2,t} =
i_2(a^2_t/E^2_{\Rang(f)})$.  It is enough to choose for each $t
\in I$ an isomorphism $g_t$ from $M^1_{i_1(a^1_2,\Dom(f))+1} \restriction
(a^1_t/E^1_{\Dom(f)})$ onto $M^2_{i_2(a^2_t,\Rang(f))+1}
\restriction (a^2_t/E^2_{\Rang(f)})$ such that: 
if $(A,B)$ is a cut of $M^1_{i_1(a^1_t,\Dom(f))}
\restriction (a_1/E^1_{\Dom(f)})$ of cofinality
$(\lambda,\lambda)$ then for some $a \in M^0_1$ we have $A < a < B$
iff for some $b \in M^2_2$ we have $g_t(A) <_{M^1_2} b <_{M_2}
g_t(B)$.  This is done as in the proof of $(*)_0$ above.]

Together it follows that $M^1_\kappa \cong M^2_\kappa$ as required.
\end{PROOF} 
\newpage

\section {Independent theories lack limit models} 

\noindent
Considering \S1 it is a natural to ask:

\begin{question}
\label{2.1}
1) Is there an unstable $T$ for which the conclusion of \ref{1.1} fails?

\noindent
2) For which unstable $T$ does the conclusion of \ref{1.1} fail?
\end{question}

\begin{remark}  1) We shall consider also relatives
$\Pr_{\lambda,\kappa}(\bar M),\Pr_{\lambda,\kappa}(T)$.

\noindent
2) In Definition \ref{2.2.7} below if 
$2^\lambda = \lambda^+$ we can restrict ourselves to $\bar M$
such that $\bigcup\{M_\alpha:\alpha < \lambda^+\} \in 
\EC_{\lambda^+}(T)$ is saturated.  The union is unique (for
$\lambda$) and there is $\bold F$ as in \ref{y.5}(3) guaranteeing this.
\end{remark}

We first note that for some $T$'s there are non-existence result (see
definitions after the claim).
\begin{theorem}
\label{2.2}  
1) If $T$ has the strong independence property
(see below, e.g. $T$ is the theory of random graphs), $|T| \le \lambda$
and $\lambda^\kappa < 2^\lambda$ \then \, $T$ does not
have a $(\lambda,\kappa)$-wk-limit model.

\noindent
2) Moreover for every $\bold F$ as in Definition \ref{y.5}(3),
there is a $\prec$-increasing continuous sequence $\bar M =
\langle M_\alpha:\alpha <
\lambda^+\rangle$ of members of $\EC_\lambda(T)$ obeying $\bold F$ such
that if $\cf(\delta_1) = \kappa = \cf(\delta_2)$
\then \, $M_{\delta_1} \cong M_{\delta_2} \Leftrightarrow \delta_1 = \delta_2$.
\end{theorem}

\begin{definition}
\label{2.2.1}  
$T$ has the strong independence
property (or is strongly independent) 
\when \,: for some $\varphi(\bar x,y) \in \bbL(\tau_T)$ for every
$M \in \EC(\tau_T)$ and pairwise distinct $a_0,\dotsc,a_{2n-1}
\in M$ for some $\bar a \in {}^{\ell g(\bar x)}M$ we have $M \models
``\varphi[\bar a,a_\ell]^{\text{if}(\ell\text{ is even})}$". 
\end{definition}

\begin{definition}
\label{2.2.7}  
Recall $S^\lambda_\kappa =: \{\delta < \lambda:\text{cf}(\delta) 
= \kappa\}$.

\noindent
1) Let $\Pr_{\lambda,\kappa}(\bar M)$ 
mean that $\bar M = \langle M_\alpha:\alpha < \lambda^+\rangle$
is $\prec$-increasing continuous, each $M_\alpha$ is of cardinality
$\lambda$ and for some club $E$ of $\lambda^+$, if $\alpha \in
S^{\lambda^+}_\lambda \cap E$ and $\delta_1 \ne \delta_2 \in
S^{\lambda^+}_\kappa \cap E$ but $\alpha < \delta_1 < \delta_2$
\then \, there is no automorphism $\pi$ of $M_\alpha$ which maps
$\{\tp(\bar a,M_\alpha,M_{\delta_1}):\bar a \in {}^{\omega >}
(M_{\delta_1})\}$ onto $\tp(\bar a,M_\alpha,M_{\delta_2}):
\bar a \in {}^{\omega >}
(M_{\delta_2})\}$ (actually even demanding just 
$\alpha \in E$ is O.K., i.e. we can prove it); note that $\pi$ acts of
$M_\alpha$ hence on $\bold S^{< \omega}(M_\alpha)$ and $\pi$ is not
necessarily the identity.

\noindent
1A) Let $\Pr_\lambda(\bar M)$ mean $\Pr_{\lambda,\lambda}(\bar M)$,
similarly for the versions below.

\noindent
2) Let $\Pr_{\lambda,\kappa}(T)$ means: for some $\bold F$ as in
\ref{y.5}(3), if $\bar M = \langle M_\alpha:\alpha <
\lambda^+\rangle$ obeys $\bold F$ then $\Pr_{\lambda,\kappa}(\bar M)$. 

\noindent
3) Let $\Pr^2_{\lambda,\kappa}(\bar M)$ be defined as in part (1) but
$\pi$ is an isomorphism from $M_{\delta_1}$ onto $M_{\delta_2}$
mapping $M_\alpha$ onto itself.  We define $\Pr^2_{\lambda,\kappa}(T)$
as in part (2) using $\Pr^2_{\lambda,\kappa}(\bar M)$.

\noindent
4) Let $\Pr^1_{\lambda,\kappa}(-)$ mean $\Pr_{\lambda,\kappa}(-)$. 
\end{definition}

\begin{remark}
\label{2.2.7t}  
1) Clearly $\Pr^2_{\lambda,\kappa}(\bar M)
\Rightarrow \Pr^1_{\lambda,\kappa}(\bar M)$ and
$\Pr^2_{\lambda,\kappa}(T) \Rightarrow \Pr^1_{\lambda,\kappa}(T)$.

\noindent
2) Also there is no point (in \ref{2.2.7}(1)) to use
$\alpha_1,\alpha_2$ as some $\bold F$ guarantee that $\alpha_1 <
\alpha_2 < \delta \in S^\lambda_\kappa$ implies there is an
automorphism of $M_\delta$ mapping $M_{\alpha_1}$ onto $M_{\alpha_2}$.
\end{remark}

\begin{PROOF}{\ref{2.2}}
  1) Assume that $\varphi(\bar x,y)$ exemplifies the strong
independence property.

For every $M \in \EC_\lambda(T)$ and function
$\bold F$ as in \ref{y.5}(3) we can find a sequence $\langle M_\alpha:\alpha <
\lambda^+\rangle$ obeying $\bold F$ such that $M \prec M_0$ and:
\mn
\begin{enumerate}
\item[$\circledast$]    if $\alpha < \lambda^+$ 
then for some $\bar c^\alpha \in {}^{\ell g(\bar x)}(M_{\alpha +1})$ we have: 
in $M_{\alpha +1}$ every $a \in M_\alpha$ satisfies
$\varphi(\bar c^\alpha,a) \Leftrightarrow a \in M$.
\end{enumerate}
\mn
Now for any $\delta < \lambda^+$ of cofinality $\kappa$ let $\langle
\alpha^\delta_\varepsilon:\varepsilon <\kappa\rangle$ be increasing with
limit $\delta$ then $\bar{\bold c}^\delta = \langle \bar
c^{\alpha^\delta_\varepsilon}:\varepsilon < \kappa\rangle$ is a
sequence of $\ell g(\bar x)$-tuples from $M_\delta$ of length
$\kappa$, and for every $a \in M_\delta$ we have:
\mn
\begin{enumerate}
\item[$(*)$]   a realizes the type $p(y,\bar{\bold c}^\delta) =
\{\varphi(\bar c^{\alpha_\varepsilon},y):\varepsilon <
\kappa\}$ in $M_\delta$ iff $a \in M$.
\end{enumerate}
\mn
The number of isomorphism types of $\tau_T$-models $M'$ of cardinality
$\lambda$ is $2^\lambda$ whereas the
number of $\langle \bar c^\alpha_i:i < \kappa \rangle$ for a given $M'$ 
is $\le \lambda^\kappa < 2^\lambda$.

For a given $\bold F$ the construction above works for every $M \in
\text{ EC}_\lambda(T)$, but $\dot I(\lambda,T) = 2^\lambda$, see \ref{y17} as
$\lambda \ge |T| + \aleph_1$ so we can finish easily, or see more in
part (2).

\noindent
2)  We can make the counterexample more explicit.  For a model $M$ and
$\bar c^\varepsilon \in {}^{\ell g(\bar x)}M$ for $\varepsilon
<\kappa$ we define $N = N[M,\langle \bar c^\varepsilon:\varepsilon <
\kappa \rangle]$ as the following submodel of $M$ (if well defined): it
is the submodel with universe the set $A = \{d \in M:M \models \varphi
[\bar c^\varepsilon,d]$ for every $\varepsilon < \kappa\}$; (note that
$N$ is not necessarily an elementary submodel of $M$ or even well
defined, e.g. $A = \emptyset$ or $A$ not closed under functions of $M$).  For
$M \in \EC_\lambda(T)$ let ${\cM}[M] = \{N \prec M:N$ is
$N[M,\langle \bar c^\varepsilon:\varepsilon < \kappa\rangle]$ for some
$\bar c^\varepsilon \in {}^{\ell g(\bar x)}M$ for $\varepsilon < \kappa\}$.
Fixing $\bold F$ as in \ref{y.5}(3) we can choose $M_\alpha \in \EC_\lambda(T)$ 
with universe $\lambda \times (1 + \alpha)$ such that
\mn
\begin{enumerate}
\item[$(*)_1$]   if $\alpha = 4 \beta +3$ and $\delta \le 4 \beta$
then $M_\alpha$ is not isomorphic to $N \prec M_\delta$ whenever there
are $\bar c^\varepsilon \in {}^{\ell g(\bar x)}(M_\delta)$ for
$\varepsilon < \kappa$ such that $N = N[M_\delta,\langle \bar
c^\varepsilon:\varepsilon < \kappa\rangle]$ 
\sn
\item[$(*)_2$]   for $\alpha < \beta < \lambda^+$ there is $\bar
c^\beta_\alpha \in {}^{(\ell g(\bar x))}(M_{\beta +1})$ such that for
every $a \in M_\beta$ we have $M_{\beta +1} \models \varphi[\bar
c^\beta_\alpha,a] \Leftrightarrow a \in M_\alpha$
\sn
\item[$(*)_3$]   the sequence $\langle M_{2 \alpha}:\alpha < \lambda^+
  \rangle$ obeys $\bold F$.
\end{enumerate}
\mn
As $\dot I(\lambda,T) = 2^\lambda$ and moreover for any theory $T_1
 \supseteq T$ of cardinality $\lambda$ we have $\dot I(\lambda,T_1,T)
 = 2^\lambda$ and for every $M \in \EC_\lambda(T)$, the number
of $N \in {\cM}[M]$ is $\le \lambda^\kappa < 2^\lambda$ we get
\mn
\begin{enumerate}
\item[$\boxtimes$]   for every appropriate $\bold F$ there is a
$\prec$-increasing continuous sequence $\langle M_\alpha:\alpha <
\lambda^+\rangle$ of models of $T$ as above such that if
$\delta_1 \ne \delta_2 < \lambda^+$ has cofinality $\kappa$ \then \,
$M_{\delta_1},M_{\delta_2}$ are not isomorphic.
\end{enumerate}
\mn
[Why?  Without loss of generality $\delta_1 < \delta_2$, let $\langle
\alpha^{\delta_2}_\varepsilon:\varepsilon < \kappa \rangle$ be increasing
with limit $\delta_2$, all $> \delta_1 +4$.  Now by $(*)_2$ we know that
$\langle \bar c^{\alpha_\varepsilon}_{\delta_1+3}:\varepsilon < \kappa \rangle$
exemplified that in $M_{\delta_2}$ there is a sequence $\langle \bar
c^\varepsilon:\varepsilon < \kappa\rangle$ which define
$M_{\delta_1+3}$, i.e. $M_{\delta_1 +3} = N[M_{\delta_2},\langle \bar
c^\varepsilon:\varepsilon < \kappa\rangle]$.

So if $M_{\delta_1} \cong M_{\delta_2}$ then there are $\bar
d^\varepsilon \in {}^{\ell g(\bar x)}(M_{\delta_1})$ for $\varepsilon
< \kappa$ such that $N[M_{\delta_1},\langle \bar
d^\varepsilon:\varepsilon < \kappa\rangle]$ is well defined and
isomorphic to $M_{\delta_1+3}$.  But consider the choice of 
$M_{\delta_1+3}$, clearly $(*)_1$ says that this is impossible.]
\end{PROOF}

\begin{observation}
\label{2.2.2}  If, inside the proof of \ref{2.2},  
in the definition of ${\cM}[M]$ we restrict
ourselves to $\langle \bar c^\varepsilon:\varepsilon < \kappa\rangle$
such that $(\forall a \in M)(\exists \varepsilon  < \kappa)(\forall
\zeta)(\varepsilon < \zeta < \kappa \rightarrow M \models
\varphi[\bar c^\varepsilon,a] \equiv \varphi[\bar c^\zeta,a])$
\then \, we can replace $\lambda^\kappa < 2^\lambda$ by
$\bold U_\kappa(\lambda) < 2^\lambda$, see \ref{y23}.
\end{observation}
\bigskip

\noindent
Considering \ref{2.2.2} (and \ref{1.1}), it is natural to ask:

\begin{question}
\label{2.2.3}  Is the independence property enough to
imply no limit models?

The problem was that the independence we can get may be
``hidden", ``camouflaged" by other ``parts" of the model.

Working harder (than in \ref{2.2}), the answer is yes.
\end{question}

\begin{theorem}
\label{2p.3}  Assume $T$ is independent.

\noindent
1) If $|T| \le \lambda = \lambda^\theta = 2^\kappa > \theta =
\cf(\theta)$ \then \, $T$ has no $(\lambda,\theta)-\md$-limit model.

\noindent
2) Moreover, there is $\bold F$ such that
\mn
\begin{enumerate}
\item[$(a)$]   $\bold F$ is a function with domain 
$\bigcup\{K_\alpha:\alpha < \lambda^+ \odd\}$ where $K_\alpha = \{M:M$ a
model of $T$ with universe $\lambda \times (1 + \alpha)\}$
\sn
\item[$(b)$]   if $\alpha < \lambda^+$ is odd and
$M \in K_\alpha$ then $M \prec \bold F(M) \in K_{\alpha +1}$
\sn
\item[$(c)$]   if $M_\alpha \in K_\alpha$ for $\alpha < \lambda^+$
is $\prec$-increasing continuous and $M_{2 \alpha +2} = \bold F(M_{2
\alpha+1})$ for $\alpha < \lambda^+$ \then \, 
for no $\alpha < \lambda^+$ is the set
$\{\delta:M_\delta \cong M_\alpha \text{ and } \cf(\delta)=\theta\}$
stationary.
\end{enumerate}
\mn
3) We can strengthen part (2) by adding in clause (c):
\mn
\begin{enumerate}
\item[$(*)$]   there are $\bar c_\alpha \in {}^\kappa(M_{2
\alpha+2})$ for $\alpha < \lambda^+$ such that: if $\langle
\alpha_{\ell,\varepsilon}:\varepsilon < \theta \rangle$ is an
increasing continuous sequence of ordinals $< \lambda^+$ with limit
$\alpha_\ell$ for $\ell=1,2$ and $\alpha_1 \ne \alpha_2$ then
there is no isomorphism $f$ from $M_{\alpha_1}$ onto $M_{\alpha_2}$
mapping $\bar c_{\alpha_{1,\varepsilon}}$ to
$\bar c_{\alpha_{2,\varepsilon}}$ for $\varepsilon < \theta$.
\end{enumerate}
\mn
4) In part (2) we can replace $K_\alpha$ (for $\alpha < \lambda^+$) by
$K_{< \lambda^+} := \cup\{K_\alpha:\alpha < \lambda^+\}$.
\end{theorem}

\begin{remark}
\label{2p.3.7}  
1) How does $2^\kappa = \lambda$ help us?

We shall consider $M_\alpha \in K_\alpha$ for $\alpha < \lambda^+$
 which is $\prec$-increasing.
We fix a sequence $\langle \bar a_t:t \in I\rangle$ in
$M_0$ such that $\langle \varphi(x,\bar a_t):t \in I\rangle$ is an
independent set of formulas (actually $I = \lambda$).  Now for any
sequence $\langle \eta_i:i < \kappa + \kappa \rangle$ of members of
${}^I 2$, and $\prec$-extension $M$ of $M_0$ we can find $N,\bar c$
such that $M \prec N,\bar c = \langle c_i:i < \kappa + \kappa\rangle$
and $N \models \varphi[c_i,\bar a_t]^{\iif(\eta_i(t))}$.  Specifically if 
$M_{2 \alpha +1}$ is already chosen then when choosing 
$M_{2 \alpha +2}$ we choose also a
sequence $\langle \eta^\alpha_i:i < \kappa + \kappa\rangle$, of
members of ${}^I 2$ and $\langle c^\alpha_i:i < \kappa +
\kappa\rangle$ such that $M_{2 \alpha +2} \models
\varphi[c^\alpha_i,\bar a_t]^{\iif(\eta^\alpha_i(t))}$.

We may look at it as coding a sequence of $\lambda$ subsets of
$\kappa$.  We essentially like to gain some information on $\langle
\eta^\alpha_i:i < \kappa + \kappa\rangle$ from $(M_{2 \alpha
+1},M_{2\alpha+2},\bar c^\alpha)$, but we are not given who are the
$\bar a_t$'s.  We shall try to use $\langle c^\alpha_i:i <
\kappa\rangle$, to distinguish between the ``true" $\bar a_t$'s and
``fakers".  We do an approximation: some will be ``exposed fakes",
which we can discard, and the others are ``perfect fakers", i.e., they
immitate perfectly some $a_t$, so it does not matter.
\end{remark}

\noindent
Clearly it suffices to prove part (3) of \ref{2p.3} for having parts
(1),(2) because $\lambda = \lambda^\kappa$ 
and the proof of part (4) is similar.
The proof is broken to some definitions and claims.
\begin{definition}
\label{2p.4}  
1) Assume $\varphi = \varphi(x,\bar y) \in \bbL(\tau_T)$ has 
the independence property in $T$. We say
$(M,\bar{\bold a})$ is a $(T,\varphi)$-candidate or an
$(I,T,\varphi)$-candidate \when \,:
\mn
\begin{enumerate}
\item[$(a)$]  $M$ is a model of $T$
\sn
\item[$(b)$]  $\bar{\bold a} = \langle \bar a_t:t \in I\rangle,
\bar a_t \in {}^{\ell g(\bar y)}M$ and $I$ is an infinite linear order  
\sn
\item[$(c)$]   $\bar{\bold a}$ is an indiscernible sequence in $M$
\sn
\item[$(d)$]  $\{\varphi(x,\bar a_t):t \in I\}$ is independent in
$M$; that is for every $\eta \in \fin(I) := \{\eta:\eta \in {}^J
2$ for some finite $J \subseteq I\}$, there is $b \in M$ such that $t \in
\Dom(\eta) \Rightarrow M \models \varphi[b,\bar a_t]^{\text{if}(\eta(t))}$.
\end{enumerate}
\mn
2) If $(M,\bar{\bold a})$ is an $(I,T,\varphi)$-candidate \then \, let
$\Gamma_{M,\bar{\bold a}} = \Gamma^{T,\varphi}_{M,\bar{\bold a}} = 
\Gamma^{T,\varphi,1}_{M,\bar{\bold a}} \cup
\Gamma^{T,\varphi,2}_{M,\bar{\bold a}}$ be the
following set of first order sentences and $\tau^+_M$ 
be the following vocabulary
\mn
\begin{enumerate}
\item[$(a)$]   $\tau^+_M = \tau_T \cup\{c:c \in M\} \cup \{{P}\}$
where $P$ is a unary predicate ($\notin \tau_T$ of course) and each 
$c \in M$ serves as an individual constant $(\notin \tau_T)$
\sn
\item[$(b)$]   $\Gamma^{T,\varphi,1}_{M,\bar{\bold a}} = \Th(M,c)_{c \in M}$
\sn
\item[$(c)$]  $\Gamma^{T,\varphi,2}_{M,\bar{\bold a}} = 
\{(\exists x)[P(x) \wedge \bigwedge\limits_{t \in J}
\varphi(x,\bar a_t)^{\iif(\eta(t))}]$: for some finite $J \subseteq
I$ and $\eta \in {}^J 2\}$ 
(so the vocabulary is $\subseteq \tau^+_M$).
\end{enumerate}
\mn
3) In (2) let $\Omega_{M,\bar{\bold a}} =
\Omega^{T,\varphi}_{M,\bar{\bold a}}$ be the family of consistent sets
$\Gamma$ of sentences in $\bbL(\tau^+_M)$ such that $\Gamma$ is of the 
form $\Gamma_{M,\bar{\bold a}}$
union with a subset of $\Phi_{M,\bar{\bold a}} = \{\neg(\exists x)[P(x)
\wedge \psi(x,\bar c) \wedge \bigwedge \limits_{t \in J} \varphi(x,\bar
a_t)^{\eta(t)}]:J \subseteq I$ is finite, $\eta \in {}^J 2,
\bar c \in {}^{\ell g(\bar z)} M$ and $\psi(x,\bar z) \in \bbL(\tau_T)\}$.

\noindent
4) For $\Gamma \in \Omega_{M,\bar{\bold a}}$ let
\mn
\begin{enumerate}
\item[$(a)$]    $\bold S_\Gamma = \{p:p \in \bold S(M)$ 
and $\Gamma \cup \{(\exists x)(P(x) \wedge \psi(x,\bar c)):
\psi(x,\bar c) \in p(x)\}$ is consistent$\}$
\sn
\item[$(b)$]   for $J \subseteq I$
and $\eta \in {}^J 2$ let $\bold S_{\Gamma,\eta} = \{p \in \bold S_\Gamma:p$
include $q^\eta_{M,\bar{\bold a}}\}$
\end{enumerate}
\mn
where 
\mn
\begin{enumerate}
\item[$(c)$]  $q^\eta_{M,\bar{\bold a}} := 
q^{T,\varphi,\eta}_{M,\bar{\bold a}} =
\{\varphi(x,\bar a_t)^{\iif(\eta(t))}:t \in \Dom(\eta)\}$.
\end{enumerate}
\mn
5) For $\Gamma \in \Omega_{M,\bar{\bold a}},\psi(x,\bar z) \in 
\bbL(\tau_T)$ and $\bar c \in {}^{\ell g(\bar z)}M$ let 

\begin{equation*}
\begin{array}{clcr}
\Xi_{M,\bar{\bold a},\Gamma}(\psi(x,\bar c)) = \{\eta \in 
\text{\rm fin}(I):&\Gamma \text{ is consistent with} \\
  &(\exists x)[P(x) \wedge \psi(x,\bar c) \wedge 
\bigwedge\limits_{t \in \Dom(\eta)} \varphi(x,\bar a_t)^{\iif(\eta(t))}]\}.
\end{array}
\end{equation*}
\end{definition}

\begin{remark}
\label{2p.6g}
1) $\fin(I) = \{\eta:\eta$ is a function from some finite $J \subseteq
   I$ to $\{0,1\}\}$.

\noindent
2) In parts (3) and (4) we could have used only $\psi(x,\bar
z) \in \{\varphi(x,\bar y),\neg \varphi(x,\bar y)\}$.
\end{remark}

\begin{observation}
\label{2p.5}  
Let $(M,\bar{\bold a})$ be a $(T,\varphi)$-candidate.

\noindent
1) $\Gamma_{M,\bar{\bold a}} \in \Omega_{M,\bar{\bold a}}$, i.e.,
$\Gamma_{M,\bar{\bold a}}$ is consistent so $\Omega_{M,\bar a}$ is non-empty.

\noindent
2) $\Omega_{M,\bar{\bold a}}$ is closed under increasing (by
$\subseteq$) unions.

\noindent
3) Any member of $\Omega_{M,\bold{\bar a}}$ can be extended to a maximal
member of $\Omega_{M,\bar{\bold a}}$.

\noindent
4) If $M \prec N$ \then \, $(N,\bar{\bold a})$ is a
$(T,\varphi)$-candidate and for every $\Gamma \in 
\Omega_{M,\bold{\bar a}}$ the set $\Gamma \cup 
\Gamma_{N,\bar{\bold a}}$ belongs to $\Omega_{N,\bar a}$.

\noindent
5) If $\langle I_\alpha:\alpha \le \delta\rangle$ is an increasing
continuous sequence of linear orders and 
$\langle N_\alpha:\alpha \le \delta\rangle$ is
$\prec$-increasing continuous sequence of models of $T,\bar{\bold a} =
\langle \bar a_t:t \in I_\delta\rangle$ and
$(N_\alpha,\bar{\bold a} \restriction I_\alpha)$ is 
a $(T,\varphi)$-candidate for $\alpha < \delta$ \then \,
$(N_\delta,\bar{\bold a})$ is a $(T,\varphi)$-candidate.

\noindent
6) In part (5), if $\Gamma_\alpha \in \Omega_{N_\alpha,\bar{\bold
a}}$ for $\alpha < \delta$ is increasing continuous with $\alpha$ 
\then \,
$\Gamma_\delta := \bigcup\{\Gamma_\alpha:\alpha < \delta\}$ belongs to
$\Omega_{N_\delta,\bar{\bold a}}$.

\noindent
7) In part (6) if $\Gamma_\alpha$ is maximal in
$\Omega_{N_\alpha,\bold a}$ for each $\alpha < \delta$ \then \,
$\Gamma_\delta$ is maximal in $\Omega_{N_\delta,\bar{\bold a}}$.

\noindent
8) If $\Gamma \in \Omega_{M,\bar{\bold a}},\psi(x,\bar z) \in 
\bbL(\tau_T),\bar c \in {}^{\ell g(\bar z)}M$ and $M \models 
(\exists x)\psi(x,\bar c)$ \then \,
\mn
\begin{enumerate}
\item[$(a)$]    the empty function belongs to 
$\Xi_{M,\bar{\bold a},\Gamma}(\psi(x,\bar c))$
\sn
\item[$(b)$]    if $I_1 \subseteq I_2$ are finite subsets of $I$ and
$\eta \in \Xi_{M,\bar{\bold a},\Gamma}(\psi(x,\bar c)) \cap
{}^{(I_1)}2$ \then \, there is $\nu \in \Xi_{M,\bar{\bold a},\Gamma}
(\psi(x,\bar c)) \cap {}^{(I_2)}2$ extending $\eta$.
\end{enumerate}
\end{observation}

\begin{PROOF}{\ref{2p.5}}
Straightforward.  
\end{PROOF}

\begin{claim}
\label{2p.6}  
Assume that $(M,\bar{\bold a})$ is a
$(T,\varphi)$-candidate and $\Gamma \in
\Omega^{T,\varphi}_{M,\bar{\bold a}}$ is maximal. 

\noindent
1) If $\psi(x,\bar z) \in \bbL(\tau_T)$ and $\bar c \in {}^{\ell
g(\bar z)}M$ and $\eta \in \Xi_{M,\bar a,\Gamma}(\psi(x,\bar c))
\subseteq \,\fin(I)$ \then \, for some $\nu$ we have $\eta \subseteq \nu \in
\fin(I)$ and $\nu \notin \Xi_{M,\bar{\bold a},\Gamma}(\neg\psi(x,\bar c))$.

\noindent
2) For every $\eta \in {}^I 2$ there are $N,b$ such that:
\mn
\begin{enumerate}
\item[$(a)$]   $M \prec N$ (and $\|N\| \le \|M\| + |T|$)
\sn
\item[$(b)$]  $b \in N$
\sn
\item[$(c)$]   if $t \in I$ then 
$N \models \varphi[b,\bar a_t]^{\iif(\eta(t))}$
\sn
\item[$(d)$]   if $\bar a \in {}^{\ell g(\bar z)}M,\psi =
\psi(x,\bar z) \in \bbL(\tau_T)$ and 
$\psi(x,\bar a) \in \,\tp(b,M,N)$
\then \, $\Gamma$ is disjoint to $\{\neg(\exists x)[P(x) \wedge
\psi(x,\bar a) \wedge \bigwedge\limits_{t \in J} \varphi(x,\bar
a_t)^{\iif(\eta(t))}]:J \subseteq I$ finite$\}$.
\end{enumerate}
\end{claim}

\begin{PROOF}{\ref{2p.6}} 
1) Assume that the conclusion fails.  Consider the
formula $\psi'(x,\bar c) := \bigwedge\limits_{t \in \Dom(\eta)}
\varphi(x,\bar a_t)^{\text{if}(\eta(t))} \rightarrow \neg \psi(x,\bar c)$.

By the assumption of the claim + the assumption toward contradiction it follows
that ``$\rho \in \fin(I) \Rightarrow \Gamma \cup \{(\exists
x)[P(x) \wedge \bigwedge\limits_{t \in \Dom(\rho)}
\varphi(x,\bar a_t)^{\iif(\rho(t))} 
\wedge \psi'(x,\bar c)]\}$ is consistent).  

\noindent
[Why?  Just note that it is enough to consider $\rho \in \fin(I)$ 
such that $\Dom(\eta) \subseteq \rho$ and we split to two
cases: \underline{first} when $\rho \restriction \Dom(\eta) \ne \eta$
then $\psi'(x,\bar c)$ adds nothing in the conjunction (and use
\ref{2p.4}(2)(c)); \underline{second} when $\rho \restriction 
\Dom(\eta) = \eta$ and we use the assumption toward the contradiction.]

So if $N'$ is a model of $\Gamma$ and we define $N''$ as $N'$ by
replacing $P^{N'}$ by $P^{N''} = \{b \in P^{N'}:N' \models
\psi'[b,\bar c]\}$ we see that $\Gamma \cup\{\neg(\exists x)[P(x)
\wedge \neg \psi'(x,\bar c)]\} \in \Omega_{M,\bar{\bold a}}$.  By the
maximality of $\Gamma$ it follows that $\neg(\exists x)[P(x) \wedge
\neg \psi'(x,\bar c)] \in \Gamma$.  But this contradicts the
assumption $\eta \in \Xi_{M,\bar{\bold a},\Gamma}(\psi(x,\bar c))$.

\noindent
2) Easy. 
\end{PROOF}

\begin{claim}
\label{2p.7}  Assume that
\mn
\begin{enumerate}
\item[$(a)$]  $(M,\bar{\bold a})$ is an $(I,T,\varphi)$-candidate
\sn
\item[$(b)$]  $\bar \eta = \langle \eta_i:i < i(*)\rangle$ and
$\eta_i \in {}^I 2$ for $i <i(*)$
\sn
\item[$(c)$]  $j(*) \le i(*)$ 
\sn
\item[$(d)$]  $\{\eta_i:i < j(*)\}$ is a dense subset of ${}^I 2$.
\end{enumerate}
\mn
\Then \, we can find $N,\bar{\bold c}$ such that
\mn
\begin{enumerate}
\item[$(\alpha)$]  $M \prec N$ and $\|N\| \le \|M\| + |T| + |i(*)|$
\sn
\item[$(\beta)$]  $\bar{\bold c} = \langle c_i:i < i(*)\rangle$ and 
$c_i \in N$
\sn
\item[$(\gamma)$]  if $i<i(*)$ and $t \in I$ then $N \models
\varphi[c_i,\bar a_t]^{\iif(\eta_i(t))}$
\sn
\item[$(\delta)$]  for every $\bar a \in {}^{\ell g(\bar y)}M$ at
least one of the following holds:
\begin{enumerate}
\item[$(i)$]  [the perfect fakers]  
for some $t \in I$ for every $\rho_0 \in \fin(I \backslash \{t\})$ 
we can find $\rho_1 \in \fin(I \backslash \{t\})$ extending $\rho_0$
such that: 
$\rho_1 \subseteq \eta_i \wedge
i < i(*) \Rightarrow N \models ``\varphi[c_i,\bar a] \equiv 
\varphi[c_i,\bar a_t]"$, i.e. for ``most" $i < i(*),\bar a,\bar a_t$
are similar
\sn
\item[$(ii)$]  [the rejected $\bar a$'s] for no $t \in I$ do 
we have $i < j(*) \Rightarrow N \models 
``\varphi[c_i,\bar a] \equiv \varphi[c_i,\bar a_t]"$.
\end{enumerate}
\end{enumerate}
\end{claim}

\begin{PROOF}{\ref{2p.7}}
By \ref{2p.3}(1), $\Gamma_{M,\bar{\bold a}} \in
\Omega_{M,\bar{\bold a}}$ hence by \ref{2p.3}(3) there is a maximal
$\Gamma \in \Omega_{M,\bar{\bold a}}$.

\noindent
Let $N,\langle c_i:i < i(*)\rangle$ be such that
\mn
\begin{enumerate}
\item[$\circledast$]  $(a) \quad M \prec N$ and $\|N\| = \|M\| +
|T| + |i(*)|$
\sn
\item[${{}}$]  $(b) \quad$ for $i < i(*),c_i \in N$ realizes some
$p_i \in \bold S_{\Gamma,\eta_i}$ (see Definition \ref{2p.4}(4)(b)).
\end{enumerate}
\mn
Clearly clauses $(\alpha),(\beta),(\gamma)$ of the desired conclusion
hold, and let us check clause $(\delta)$.  So assume that $\bar a \in
{}^{\ell g(\bar y)}M$ and clause (ii) there fails so we can choose $t
\in I$ such that $i < j(*) \Rightarrow N \models ``\varphi[c_i,\bar a]
\equiv \varphi[c_i,\bar a_t]"$.

So it is enough to prove clause (i) for $t$; 
toward this assume $\rho_0 \in \fin(I)$ satisfies 
$t \notin \Dom(\rho_0)$, i.e. $\rho_0 \in \fin(I \backslash \{t\})$.  
Let $\rho_1 \in \fin(I)$ extend $\rho_0$ be such that $\rho_1(t) = 0$.
By clause (d) of the assumption we know that for some $i < j(*)$ we have
$\rho_1 \subseteq \eta_i$ but (see above) $N \models
``\varphi[c_i,\bar a] \equiv \varphi[c_i,\bar a_t]"$
 hence $\rho_1 \in \Xi_{M,\bar{\bold
a},\Gamma}(\varphi(x,\bar a)^{\iif(\eta_i(t))})$ which means that $\rho_1
\in \Xi_{M,\bar a,\Gamma}(\varphi(x,\bar a)^{\text{if}(\rho_1(t))}$).  
Now apply claim \ref{2p.6}(1)
to $\psi(x,\bar c) := \varphi(x,\bar a)^{\iif(\rho_1(t))}$, so we know that
for some $\nu$ we have $\rho_1 \subseteq \nu \in \fin(I)$ and
$\nu \notin \Xi_{M,\bar{\bold a},\Gamma}(\neg \varphi(x,a)^{\iif(\rho_1(t))})$
hence 
\mn
\begin{enumerate}
\item[$(*)_1$]   if $i < i(*)$ satisfies 
$\nu \subseteq \eta_i$ then $N \models
``\varphi[c_i,\bar a]^{\iif(\rho_1(t))}"$ which means 
$N \models ``\varphi[c_i,\bar a] \equiv \varphi[c_i,\bar a_t]"$.
\end{enumerate}
\mn
Let $\rho_2 \in {}^{\Dom(\nu)}2$ be such that $\rho_2(t)=1$ and
$s \in \Dom(\nu) \backslash \{t\} \Rightarrow \rho_2(s) =
\nu(s)$.  We repeat the use of \ref{2p.6}(1) for $\rho_2$ instead of
$\rho_1$ and get $\nu'$ such that $\rho_2 \subseteq \nu' \in \fin(I)$ and
\mn
\begin{enumerate}
\item[$(*)_2$]   if $i < i(*)$ satisfies 
$\nu' \subseteq \eta_i$ then $N \models ``\neg
\varphi[c_i,\bar a]^{\iif(\rho_2(t))}"$ which means that $N \models ``\varphi
[c_i,\bar a] \equiv \varphi[c_i,\bar a_t]"$.
\end{enumerate}
\mn
Let $\rho_3 = \nu' \restriction (\Dom(\nu') \backslash \{t\})$
and by $(*)_1 + (*)_2$ the function $\rho_3 \in \fin(I)$ is as
required in subclause $(i)$ (for our $\bar a,t,\rho_0$) in clause
$(\delta)$ of the claim, so we are done.  
\end{PROOF}

\begin{definition}
\label{2p.11}  
1) For a model $M$ of $T$, formula $\varphi =
\varphi(x,\bar y) \in \bbL(\tau_M),\bar{\bold c} = \langle
c^\varepsilon_i:i < i(*),\varepsilon < \theta\rangle$ such that
$c^\varepsilon_i \in M$ let ${\cP}_\varphi(\bar{\bold c},M) =
\{{\cU} \subseteq i(*)$: for some $\bar a \in {}^{\ell g(\bar y)}M$
for every $\varepsilon < \theta$ large enough ${\cU} = u_\varphi(\bar
a,\langle c^\varepsilon_i:i < i(*)\rangle,M)\}$ where $u_\varphi(\bar
a,\langle c^\varepsilon_i:i < i(*)\rangle,M) = \{i < i(*):M 
\models \varphi[c^\varepsilon_i,\bar a]\}$.

\noindent
2) For a model $M$ and $\varphi = \varphi(x,\bar y) \in 
\bbL(\tau_M)$ let ${\bold P}^\varphi_{i(*),\theta}(M) =
\{{\cP}_\varphi(\bar{\bold c},M):\bar{\bold c}$ has the form $\langle
c^\varepsilon_i:i < i(*),\varepsilon < \theta\rangle$ with
$c^\varepsilon_i \in M\}$.

\noindent
3) For $M_1 \prec M_2$ and $\varphi = \varphi(x,\bar y) \in 
\bbL(\tau_{M_\ell})$ and $\bar c = \langle c_i:i < i(*)\rangle \in
{}^{i(*)}(M_2)$ let ${\cP}_\varphi(\bar{\bold c},M_1,M_2) = \{u_\varphi(\bar a,
\bar{\bold c},M_2):\bar a \in {}^{\ell g(\bar y)} M_1\}$.
\end{definition}

\begin{observation}
\label{2p.12}  
For $M,\varphi(x,\bar y),i(*),
\theta$ as in Definition \ref{2p.11}, 
$\bold P^\varphi_{i(*),\theta}(M)$ has cardinality $\le \|M\|^{|i(*)|+\theta}$.
\end{observation}

\begin{claim}
\label{2p.13}  If $I$ is a linear order of cardinality
$\le 2^\kappa$ \then \, we can find a uniform filter $D$ on $\kappa$ and a
sequence $\bar{\cU}^* = \langle {\cU}^*_{t,\ell}:t \in I,
\ell \in \{0,1,2\}\rangle$ of members of $[\kappa]^\kappa$ such that:
\mn
\begin{enumerate}
\item[$\boxtimes_1$]  $(a) \quad$ for each $t \in I,\langle{\cU}^*_{t,\ell}:
\ell = 0,1,2\rangle$ is a partition of $\kappa$
\sn
\item[${{}}$]  $(b) \quad {\cU}^*_{t,2} \in D$ for $t \in I$
\sn
\item[${{}}$]  $(c) \quad \cP(\kappa)/D$ has cardinality $2^\kappa$,
  moreover extend some free Boolean 

\hskip25pt  Algebra of cardinality $2^\kappa$
\sn
\item[$\boxtimes_2$]   if $(M,\bar{\bold a})$ is a
$(I,T,\varphi)$-candidate and $\bar{\cU} = \langle {\cU}_t:t \in
I\rangle$ satisfies $\cU^*_{t,1} \subseteq \cU_t 
\subseteq \cU^*_{t,1} \cup \cU^*_{t,2}$ \then \, for
some $(N,\bar{\bold c})$ we have
\begin{enumerate}
\item[$(a)$]  $M \prec N$ and $\|N\| \le \|M\| + |\tau_T| + \kappa$
\sn
\item[$(b)$]  $\bar{\bold c} \in {}^\kappa N$
\sn
\item[$(c)$]  if $t \in I,\bar a \in {}^{\ell g(\bar y)} M$ and
${\cU}^*_{t,1} \subseteq u(\bar a,\bar{\bold c},N) \subseteq 
{\cU}^*_{t,1} \cup {\cU}^*_{t,2}$ then 
$u(\bar a,\bar{\bold c},N) = {\cU}_t \mod D$
\sn
\item[$(d)$]  if $t \in I$ then ${\cU}_t \in \cP_\varphi(\bar{\bold c},M,N)$.
\end{enumerate}
\end{enumerate}
\end{claim}

\begin{PROOF}{\ref{2p.13}}
We replace $\kappa$ by $\kappa + \kappa$.

Let $\langle \eta^*_i:i < \kappa\rangle$ be a sequence of members
of ${}^I \kappa$ which is dense possible by \cite{EK}.

For $\ell=0,1,2$ let ${\cU}_{t,\ell} = \{i < \kappa:\eta_i(t) = \ell$
or $\eta_i(t) \ge 3 \wedge \ell=2\}$. 
Notice that it is important that $D$ is defined
independently of ${\cU}_t$ and we should therefore define it
here.  But for clarity of exposition we will only define it later.

Let (where $\alpha + \cU = \{\alpha + \beta:\beta \in \cU\}$)
\mn
\begin{enumerate}
\item[$(*)_1$]  ${\cU}^*_{t,0} = {\cU}_{t,0} \cup (\kappa + {\cU}_{t,0})$
\sn
\item[$(*)_2$]  ${\cU}^*_{t,1} = {\cU}_{t,1} \cup {\cU}_{t,2} \cup 
(\kappa + {\cU}_{t,1})$
\sn
\item[$(*)_3$]  ${\cU}^*_{t,2} = \kappa + {\cU}_{t,2}$.
\end{enumerate}
\mn
Assume $\bar{\cU} = \langle {\cU}_t:t \in I\rangle$ is such that
\mn
\begin{enumerate}
\item[$(*)_4$]  ${\cU}^*_{t,1} \subseteq {\cU}_t \subseteq {\cU}^*_{t,1} \cup
{\cU}^*_{t,2} \subseteq \kappa + \kappa$.
\end{enumerate}
\mn
Define $\eta_i = \eta^{\bar{\cU}}_i \in {}^{\kappa + \kappa}2$ for
$i < \kappa + \kappa$ by:
\mn
\begin{enumerate}
\item[$(*)_5$]  $\eta_i(t) = \begin{cases} 0 &i \notin {\cU}_t \\
  1 &i \in {\cU}_t
\end{cases}$
\end{enumerate}
\mn
Let $\bar \eta = \langle \eta_i:i < \kappa + \kappa\rangle$.

\noindent
Notice that $\langle \eta_i:i < \kappa\rangle$ is dense in ${}^I 2$ by
the choice of $\eta_i$ in $(*)_5$ and $(*)_2$ because ${\cU}_t \cap \kappa =
{\cU}^*_{t,1} \cup \cU^*_{t,2}$ and $\langle
\eta^*_i:i < \kappa\rangle$ was dense in ${}^I \kappa$.
By \ref{2p.7} applied to $(M,\bar{\bold a},\bar \eta),i(*) = \kappa
+ \kappa,j(*) = \kappa$ we can find $N,\bar{\bold c}$ as there and we
should check that they are as required.  Clauses
$(\alpha),(\beta),(\gamma)$ of the conclusion of \ref{2p.7} give
the ``soft" demands.

More specifically clause (a) of $\boxtimes_1$ holds by the choice of
the $\cU^*_{t,\ell}$'s; clauses (b),(c) of $\boxtimes_2$ holds by the
conclusion of \ref{2p.7}.

Clearly
\mn
\begin{enumerate}
\item[$(*)_6$]  $u_\varphi(\bar a_t,\bar{\bold c},N) = \{i < \kappa +
\kappa:N \models ``\varphi[c_i,\bar a_t]"\} 
= \{i < \kappa + \kappa:\eta_i(t)=1\} = {\cU}_t$
\end{enumerate}
\mn
hence
\mn
\begin{enumerate}
\item[$(*)_7$]  $t \in I \Rightarrow {\cU}_t \in 
\cP_\varphi(\bar{\bold c},M,N)$.
\end{enumerate}
\mn
So we see that demand (d) of $\boxtimes_2$ is satisfied - all the 
${\cU}_t$ are included.  We still need to prove clause (c) of 
$\boxtimes_2$, that is to show that there are no ``fakers" and, of
course, to define $D$. 

So assume
\mn
\begin{enumerate}
\item[$\odot_1$]   ${\cU}^*_{t_1} \subseteq u_\varphi(\bar a,\bar c,N) 
\subseteq {\cU}^*_{t_1} \cup {\cU}^*_{t_1}$ for some
$t_1 \in I$ and $\bar a \in {}^{\ell g(\bar y)}M$.
\end{enumerate}
\mn
Denote ${\cU} = u_\varphi(\bar a,\bar c,N)$.  We need to show ${\cU} =
{\cU}_{t_1} \mod D$.  

By clause $(\delta)$ of the conclusion
of \ref{2p.7} for $\bar a$ one of the two clauses there (i),(ii) occurs.

Recall that
\mn
\begin{enumerate}
\item[$\odot_2$]  ${\cU}^*_{t_1,1} \subseteq {\cU} \subseteq
{\cU}^*_{t_1,1} \cup {\cU}^*_{t_1,2}$.
\end{enumerate}
\mn
So ${\cU} \cap \kappa = {\cU}^*_{t_1,1} \cap \kappa = 
{\cU}_{t_1,1} \cup {\cU}_{t_1,2}$.

Now
\mn
\begin{enumerate}
\item[$\odot_3$]  for $\bar a$ clause (ii) of \ref{2p.7}$(\delta)$ fails.
\end{enumerate}
\mn
[Why?  Because $t_1$ witnesses this by the above equality and for 
each $i < \kappa$

\[
i \in {\cU} \Leftrightarrow i \in {\cU}_{t_1,1} \cup 
{\cU}_{t_1,2} \Leftrightarrow \eta_i[t_1]=1 \Leftrightarrow N \models
``\varphi[c_i,\bar a_{t_1}]".]
\]
\mn
By \ref{2p.7}$(\delta)$ and $\odot_3$ we can deduce: 
\mn
\begin{enumerate}
\item[$\odot_4$]  for $\bar a$, clause (i) of \ref{2p.7}$(\delta)$ 
holds so there is $t_2$ witnessing it.
\end{enumerate}

Next
\mn
\begin{enumerate}
\item[$\odot_5$]  $t_1 = t_2$.
\end{enumerate}
\mn
Why?  Toward contradiction assume $t_1 \ne t_2$ hence we can find $\rho_1
\in \fin(I \backslash \{t_2\})$ such that
\mn
\begin{enumerate}
\item[$\circledast_{5.1}$]  $\rho_1 \subseteq \eta_i \wedge i < \kappa
+ \kappa \Rightarrow N \models ``\varphi[c_i,\bar a] \equiv
\varphi[c_i,\bar a_{t_2}]"$.
\end{enumerate}
\mn
Without loss of generality $t_1 \in \Dom(\rho_1)$ and define
$\rho_2 = \rho \cup \{(t_2,1-\rho_1(t_1))\}$, so $\rho_1 \subseteq
\rho_2 \in \fin(I)$.  As $\{\eta_i:i < \kappa\}$ was chosen
as a dense subset of ${}^{\{0,1,2\}}I$, there is $i < \kappa$ such that $\rho_2
\subseteq \eta^*_i$, hence by $\circledast_{5.1}$
\mn
\begin{enumerate}
\item[$\circledast_{5.2}$]  $N \models ``\varphi[c_i,\bar a] \equiv
\varphi[c_i,\bar a_{t_2}]"$
\end{enumerate}
\mn
but by the choice of $\eta_i$ we have:
\mn
\begin{enumerate}
\item[$\circledast_{5.3}$]  $N \models 
``\varphi[c_i,\bar a_{t_2}]^{\iif(\eta_i(t_2))}"$
\end{enumerate}
\mn
but $\eta_i(t_2) = 1 - \rho_1(t_1)$ hence together
\mn
\begin{enumerate}
\item[$\circledast_{5.4}$]  $N \models 
``\varphi[c_i,\bar a]^{\iif(1-\rho_1(t_1))}"$
\end{enumerate}
\mn
but by the choice of $c_i$ we have:
\mn
\begin{enumerate}
\item[$\circledast_{5.5}$]  $N \models 
``\varphi[c_i,\bar a_{t_1}]^{\iif(\eta_i(t_1))}"$ 
\end{enumerate}
\mn
hence by $\circledast_{5.1}$
\mn
\begin{enumerate}
\item[$\circledast_{5.6}$]  $N \models 
``\varphi[c_i,\bar a_{t_1}]^{\iif(\rho_1(t_1))}"$.
\end{enumerate}
\mn
But $\circledast_{5.5} + \circledast_{5.6}$ contradict the choice of $t_1$ as
$i <\kappa$ using $\odot_2$ so $\odot_5$ holds, i.e. $t_1 =t_2$.]

Now subclause (i) of \ref{2p.7}$(\delta)$ tells us
\mn
\begin{enumerate}
\item[$\odot_6$]   for every $\rho_0 \in \fin(I \backslash
\{t_1\})$ there is $\rho_1 \in \fin(I \backslash \{t_1\})$
extending $\rho_0$ such that
\begin{enumerate}
\item[$(a)$]    $\rho_1 \subseteq \eta_i \wedge i < \kappa +
\kappa \Rightarrow N \models \varphi[\bar c_i,\bar a] \equiv
\varphi[c_i,\bar a_{t_1}]$ hence
\sn
\item[$(b)$]   if $\rho_1 \subseteq \eta_i \wedge i < \kappa +
\kappa \Rightarrow i \in u_\varphi(\bar a,\bar{\bold c},N)
\Leftrightarrow i \in {\cU}_{t_1}$.
\end{enumerate}
\end{enumerate}
\mn
So let

\begin{equation*}
\begin{array}{clcr}
D = \{{\cU} \subseteq \kappa + \kappa:&\text{ for every } \rho_0
\in \fin(I) \text{ there is } \rho_1, \\
  &\rho_0 \subseteq \rho_1 \in \fin(I) \text{ such that} \\
  &\kappa \le i < \kappa + \kappa \wedge \rho_1 \subseteq \eta_i
  \Rightarrow i \in {\cU}\}.
\end{array}
\end{equation*}

\mn
Clearly the filter $D$ satisfies clause $\boxtimes_1(c)$ so 
we are done.  
\end{PROOF}

\begin{PROOF}{\ref{2p.3}}
\underline{Proof of the Theorem \ref{2p.3}(3)}  
Like the proof \ref{2.2} of the case ``$T$ has the 
strong independence property."
\end{PROOF}

\begin{remark}
\label{2p.8.3}  
1) The $\bold F$ we construct works for all
$\theta = \text{ cf}(\theta) < \lambda$ for which $\lambda =
\lambda^\theta$ simultaneously.
\end{remark}

\begin{discussion}
\label{2p.9}  
Can we prove \ref{2p.3}
also for $\lambda$ strongly inaccessible?  Toward this
\mn
\begin{enumerate}
\item[$(a)$]   we have to use $\bar{\bold c}_\alpha = \langle
c_{\alpha,i}:i < \lambda\rangle$, instead $\langle c_{\alpha,i}:i
< \kappa \rangle$
\sn
\item[$(b)$]   each $M_\alpha$ has a presentation $\langle
M_{\alpha,\zeta}:\zeta < \lambda\rangle$
\sn
\item[$(c)$]    for a club $E$ of $\mu < \lambda$, we use $\langle
c_{\alpha,i}:i < \mu\rangle \char 94 \langle c_\mu\rangle$ to code
${\cU}_\alpha \cap \mu$
\sn
\item[$(d)$]   instead $i,\kappa +i$ we use $2i,2i+1$.
\end{enumerate}
\mn
So the problem is: arriving to $\mu$, we have already committed
ourselves for the coding of ${\cU}_\alpha \cap \mu'$ for $\mu' \in
E_\alpha \cap \mu$, what freedom do we have in $\mu$?

Essentially we have a set $\Lambda_\mu \subseteq {}^{2^\mu}2$
quite independent, and for $\mu_1 < \mu_2$, there is a natural
reflection, the set of possibilities in ${}^\lambda 2$ is decreasing.
But the amount of freedom left should be enough to code.  We shall
deal in \cite{Sh:906} with the inaccessible case.
\end{discussion}

\begin{question}  Can we improve \ref{2p.3}(3) in the case of $T$ not
strongly dependent? 
\end{question}

\begin{claim}
\label{2.2.14}  1) Assume $T$ has the strong
independence property.  If $\lambda \ge \kappa = \cf(\kappa),
2^{\min\{2^\kappa,\lambda\}} > \lambda^\kappa$
and $\lambda > |T|$, \then \, $\Pr_{\lambda,\kappa}(T)$.

\noindent
2) Assume $T$ is independent.  If $\lambda,\kappa$ are as above,
\then \, $\Pr(\lambda,\kappa)$. 
\end{claim}

\begin{PROOF}{\ref{2.2.14}}
1) Let $\varphi(\bar x,y)$ exemplify ``$T$ has the
strong independent property", see Definition \ref{2.2.1}.

We choose $\bold F$ such that:
\mn
\begin{enumerate}
\item[$(*)$]   if $\bold F(\langle M_i:i \le \alpha +1\rangle) \prec
M_{\alpha +2}$ then for every $i \le \alpha$ for some $\bar c = \bar
c_{\alpha,i} \in {}^{\ell g(\bar x)}(M_{\alpha +2})$ the set $\{a \in
M_i:M_{\alpha +1} \models \varphi[\bar c,a]\}$ does not belong to $\{\{a
\in M_i:M_{\alpha +1} \models \varphi[\bar d,a]\}:\bar d \in {}^{\ell
g(\bar x)}(M_{\alpha +1})\}$.
\end{enumerate}
\mn
We continue as in the proof of \ref{2.2}.

\noindent
2) Similarly (recalling the proof of \ref{2p.3}).
\end{PROOF}
\newpage

\section {More on $(\lambda,\kappa)$-limit for $T_{\ord}$} 

It is natural to hope that a $(\lambda,\kappa)$-i.md.-limit model 
is $(\lambda,\kappa)$-superlimit but in Theorem \ref{nl.11}.
we prove that there is no $(\lambda,\kappa)$-superlimit model for
$T_{\rd}$, see Definition \ref{0.2}(2). 

We conclude by showing that the $(\lambda,\kappa)$-i.md.-limit
model has properties in the direction of superlimit.  By \ref{sl.21}
it is $(\lambda,S)$-limit$^+$, that is if 
$\langle M_\alpha:\alpha < \lambda^+\rangle$ is a
$\subseteq$-increasing sequence of $(\lambda,\kappa)$-i.md-limit
models for a club of $\delta < \lambda^+$ of cofinality
$\kappa$ the model $\cup\{M_i:i < \delta\}$ is a $(\lambda,\kappa)$-i.md.-limit
model.  Also in \S1 the function $\bold F$ does not need memory.

\begin{hypothesis}
\label{nl.0}  
1) $\lambda = \lambda^{< \lambda} > \kappa = \cf(\kappa)$.

\noindent
2) We deal with $\EC_{T_{\rd}}(\lambda)$, ordered by
$\subseteq$, so $M,N$ denotes members of $\EC_\lambda(T_{\rd})$.
\end{hypothesis}

\noindent
Recall $T_{\rd}$ is from Definition \ref{0.2}(2) and recalling
Definition \ref{z15}.
\begin{definition}
\label{nl.0.1}  
1) If $M \subseteq N$ and $(C_1,C_2)$ is a cut of $M$ let
$N^{[(C_1,C_2)]} = N\restriction \{a \in N:a$ realizes the cut
$(C_1,C_2)$ of $M$ which means $c_1 \in C_1 \Rightarrow c_1 <_N a$
and $c_2 \in C_2 \Rightarrow a <_N c_2\}$.

\noindent
2) For a cut $(C_1,C_2)$ of $M,A$ is unbounded in the cut if $A \cap
   C_1$ is unbounded in $C_1$ and $A \cap C_2$ is unbounded from below
   in $C_2$.

\noindent
3) Let cut$_\kappa(M) = \{(C_1,C_2):(C_1,C_2)$ a cut of $M$ such that
   cf$(C_1,C_2) = (\kappa,\kappa)\}$ for any $M \in \EC_\lambda(T_{\rd})$.
\end{definition}
\bigskip

\noindent
\centerline{$* \qquad * \qquad *$}

\begin{definition}
\label{sl.1} 
1) We say $\bar M$ is a $(\lambda,\kappa)$-sequence \when \,:
\mn
\begin{enumerate}
\item[$(a)$]  $\bar M = \langle M_i:i \le \kappa\rangle$ is
$\subseteq$-increasing continuous sequence of members of
$\EC_\lambda(T_{\rd})$ 
\sn
\item[$(b)$]   if $i < \kappa$ and $(C_1,C_2)$ is a cut of $M_i$
then $(\alpha)$ or $(\beta)$ hold but not both where
\sn
\begin{enumerate}
\item[$(\alpha)$]  $\cf(C_1,C_2) = (\lambda,\lambda)$ and no $a
\in M_\kappa \backslash M_i$ realizes $(C_1,C_2)$
\sn
\item[$(\beta)$]   $M_\kappa \restriction \{a \in M_i:a$ realizes
$(C_1,C_2)\}$ is infinite, moreover has neither first nor last member
\end{enumerate}
\item[$(c)$]   for every $a <_{M_i} b$ the model 
$(M_\kappa \restriction (a,b)_{M_\kappa})$ is universal 
(for $\EC_\lambda(T_{\rd})$, usual embedding).
\end{enumerate}
\end{definition}

\begin{remark}  Compared to \S1 we do not require
\mn
\begin{enumerate}
\item[$(d)$]  if $i < \kappa$ and $(C_1,C_2)$ is a cut of $M_i$ not
realized by any $a \in M_\kappa$ \then \,: for every $j<i$, either $M_j$
is unbounded in $(C_1,C_2)$, or for some $a_1 \in C_1,a_2 \in
C_2$ the interval $(a_1,a_2)_{M_i}$ is disjoint to $M_j$.
\end{enumerate}
\end{remark}

\begin{claim}
\label{sl.2} 
1) If $M = \langle M_i:i \le \kappa\rangle$ is a 
$(\lambda,\kappa)$-sequence \then \,
$M_\kappa$ is $(\lambda,\kappa)$-i.md.-limit (and so for some
$\bar M' = \langle M'_i:i \le \kappa\rangle$ the statement
$\circledast_{\bar M'}$ from the proof of \ref{1.1} holds and
$M_\kappa \cong M'_\kappa$).

\noindent
2) If $(C_1,C_2)$ is a cut of $M_i,i < \kappa$ and 
$(b)(\beta)$ of Definition \ref{sl.1} holds, \then \, for some
$j \in (i,\kappa),M_j \restriction \{a \in M_j:a$ realizes 
$(C_1,C_2)\}$ is a universal model of $T_{\rd}$.

\noindent
3) If $M$ is $(\lambda,\kappa)$-i.md.-limit, \then \, there is a
   $(\lambda,\kappa)$-sequence $\langle M_i:i \le \kappa\rangle$ such
   that $M_\kappa =M$.

\noindent
4) If $M \in \EC_\lambda(T_{\rd})$ then:
\mn
\begin{enumerate}
\item[$(a)$]  if $\lambda = \|M\| = \lambda^{< \lambda}$ \then \, the
number of cuts of $M$ of cofinality $\ne (\lambda,\lambda)$ is at most
$\lambda$ 
\sn
\item[$(b)$]   if $\lambda = \|M\| = \|M\|^\kappa$ \then \, the
number of cuts of $M$ of cofinality $(\kappa,\kappa)$
is at most $\lambda$
\sn
\item[$(c)$]  if $\lambda = \|M\|$ \then \, the number of cuts of 
cofinality $(\sigma_1,\sigma_2)$ where $\sigma_1 \ne \sigma_2$ is $\le
\lambda$. 
\end{enumerate}
\end{claim}

\begin{PROOF}{\ref{sl.2}}
  1) As in the proof of \ref{1.1}, using parts (2),(3) 
see \ref{sl.2.3}(1).

\noindent
2) There are $j \in (i,\kappa)$ and 
$c \in M^{[c_1,c_2]}_j$ and $d \in M_j$ such that
$c < d,(c,d)_{M_j} \cap M_i = \emptyset$.  Now use \ref{sl.2.3} below.

\noindent
3) Should be clear.

\noindent
4) Clauses (a),(b) are easy and clause (c) holds by \cite[VIII,\S0]{Sh:c}.
\end{PROOF}

\begin{remark}
A difference between Definition \ref{sl.1} and the
earlier one is that we do not ask that a dense set of cuts of
cofinality $(\lambda,\lambda)$ of $M_i$ is realized in $\bigcup\{M_j:j < i\}$. 
\end{remark}

\begin{observation}
\label{sl.2.3}  
1) If $M \in \EC_\lambda(T_{\rd})$ is universal, 
$\lambda = \lambda^\kappa$ and $M = \bigcup\limits_{i < \kappa}
I_i$ \then \, for at least one $i < \kappa,M \restriction I_i$ is
universal for $\EC_\lambda(T_{\rd})$.

\noindent
2) If $M$ is $(\lambda,\kappa)$-i.md.-limit or just weakly
   $(\lambda,\kappa)$-i.md-limit and $a <_M b$ \then \, for some $N$:
\mn
\begin{enumerate}
\item[$(a)$]  $N \subseteq M \restriction (a,b)_M$
\sn
\item[$(b)$]  $N \in \EC_\lambda(T_{\rd})$ is universal
\sn
\item[$(c)$]   every $(C_1,C_2) \in \text{ cut}_\kappa(N)$ is
realized in $M$, (but not used).
\end{enumerate}
\end{observation}

\begin{PROOF}{\ref{sl.2.3}}
1) Let $N = {}^\kappa M$ ordered lexicographically, so $N \in
\EC_\lambda(T_{\rd})$ hence there is an embedding $f$ of $N$ 
into $M$.  We try to choose $\nu_i \in {}^i |M|$ by induction 
on $i < \kappa$ such that $j < i \Rightarrow
\nu_j \triangleleft \nu_i$ and $\nu_i \triangleleft \eta \in
{}^\kappa M \Rightarrow f(\eta) \notin I_i$ and for $i=0$ or $i$
   limit there is no problem to choose $\nu_i$.  We cannot succeed as
then $f(\bigcup\limits_{i < \kappa} \nu_i) \in M \backslash
   \bigcup_{j<i} I_j$, contradiction.  So for some $i < \kappa,\nu_i$
   has been chosen but we cannot choose $\nu_{i+1}$.  So for each $a
   \in M$ there is $\eta_a \in {}^\kappa M$ such that $\nu_i \char
   94 \langle a \rangle \triangleleft \eta_a \wedge f(\eta_a) \in
   I_i$.  So $a \mapsto f(\eta_a)$ is an embedding of $M$ into $I_i$,
   so we are done.

Alternatively, let $N \subseteq M$ be a saturated model of
$T_{\ord}$.  Try to choose $c_i <_N d_i$ by induction on $i <\kappa$
such that $j < i \Rightarrow c_j <_N < c_i <_N d_i <_N d_j$ and
$(c_{i+1},d_{i+1})_\mu \cap I_i = \emptyset$.  For some $i$ we have
$(c_i,d_i)$ well defined but we cannot choose $(c_{i+1},d_{i+1})$
hence $I_i \cap (c_i,d_i)_N$ is dense in $(c_i,d_i)_N$.

\noindent
2) Should be clear.  
\end{PROOF}

\begin{claim}
\label{sl.4.21}  If $S \subseteq S^{\lambda^+}_\kappa$
is stationary and $M \in \EC_\lambda(T_{\rd})$ is
$(\lambda,S)$-wk-limit \then \, $M$ is $(\lambda,\kappa)$-i.md.-limit.
\end{claim}

\begin{PROOF}{\ref{sl.4.21}}  
Let $\bold F_1$ witness that $M$ is 
$(\lambda,S)$-wk-limit.  We can find $\bar M = \langle M_\alpha:\alpha
< \lambda^+\rangle$ so $M_\alpha \in \EC_\lambda(T_{\oor})$
is a $\subseteq$-increasing continuous sequence
such that $\bar M$ obeys $\bold F_1$, such that in addition the sequence is as
in the proof of \ref{1.1}.  So by the choice of the set
$S' = \{\delta \in S:M_\delta \cong M\}$ is stationary, and by
\ref{1.1} the set $S'' = \{\delta:M_\delta$ is
$(\lambda,\kappa)$-i.md.-limit$\}$ is $\equiv S^{\lambda^+}_\kappa
\mod \cD_{\lambda^+}$.  Together $S' \cap S'' \ne \emptyset$
hence $M$ is $(\lambda,\kappa)$-i.md.-limit.
\end{PROOF}

\begin{definition}
\label{nl.5.7}  
1) We say that $\bar M$ witnesses that $M$ is 
$(\lambda,\kappa)$-i.md.-limit \when \,:
\mn
\begin{enumerate}
\item[$(*)$]   $\bar M = \langle M_\alpha:\alpha \le \kappa\rangle$
is such that $\circledast_{\bar M}$ from the proof of \ref{1.1}
holds and $M = M_\kappa$.
\end{enumerate}
\end{definition}

\begin{claim}
\label{nl.11}  
For $\lambda = \lambda^{< \lambda} >
\kappa = \cf(\kappa)$ \then \, there is no
$(\lambda,\kappa)$-superlimit model of $T_{\ord}$.
\end{claim}

\begin{remark}
It is trivial to show that there is no superlimit $M \in
\EC_\lambda(T)$, but we deal with $(\lambda,\kappa)$-superlimit.
\end{remark}

\begin{PROOF}{\ref{nl.11}}
Assume there is one, then by \S1 it is a
$(\lambda,\kappa)$-i.md.-limit model so there is $\bar M = \langle M_i:i \le
\kappa\rangle$ which witnesses this (i.e. such that
$\circledast^\kappa_{\bar M}$ from the proof of \ref{1.1}) hence each
$M_{i+1}$ is saturated.
As $M_0$ is universal for $\EC_\lambda(T_{\rd})$, we can find $c_\eta \in
M_0$ for $\eta \in {}^{\kappa \ge}(\lambda +1)$ such that $\eta
<_{\lex} \nu \Rightarrow c_\eta <_{M_0} c_\nu$.  For $\zeta <
\kappa$ let $\Lambda_\zeta = \{\eta \in {}^\kappa(\lambda +1)$: for
every $\varepsilon \in [\zeta,\kappa)$ we have
$\eta(\varepsilon)=\lambda\}$ and let $\Lambda_\kappa = \Lambda =
\bigcup\{\Lambda_\zeta:\zeta < \kappa\}$ so $\langle \Lambda_\zeta:\zeta <
\kappa\rangle$ is $\subseteq$-increasing.  For $\eta \in \Lambda_\kappa$ let
$(C_{1,\eta},C_{2,\eta})$ be the cut of $M_\kappa$ with $C_{1,\eta}
= \{a \in M_\kappa:a <_{M_\kappa} c_{\eta \restriction i}$ for some $i <
\kappa\}$.  So $\cf(C_{1,\eta},C_{2,\eta}) = (\kappa,\kappa)$ recalling
clause $(i)_1$ of $\circledast^\kappa_{\bar M}$ from the proof of \ref{1.1}.

Let $\langle d_j:j < \lambda\rangle$ be a decreasing sequence in $M_0$
and let
\mn
\begin{enumerate}
\item[$\circledast_0$]   $M'_i = M_i \restriction 
\{d:d_j <_M d$ for some $j < \lambda\}$ for $i \le \kappa$.
\end{enumerate}
We can choose $M_*$ such that:
\mn
\begin{enumerate}
\item[$\circledast_1$]  $(a) \quad M_\kappa \subseteq M_* 
\in \EC_{T_{\rd}}(\lambda)$
\sn
\item[${{}}$]   $(b) \quad$ if $c \in M_* \backslash M_\kappa$ then
some $\eta \in \Lambda_\kappa,c$ realizes the cut
$(C_{1,\eta},C_{2,\eta})$
\sn
\item[${{}}$]  $(c) \quad$ for every $\eta \in \Lambda$ there is an
isomorphism $f_\eta$ from $M'_\kappa$ onto
$M^{[(C_{1,\eta},C_{2,\eta})]}_*$
\sn
\item[$\circledast_2$]   for $\zeta \le \kappa$ let
$M^*_\zeta = M_* \restriction \{c:c \in M_\kappa$ or $c \in M_*$
realizes the cut 

\hskip75pt $(C_{1,\eta},C_{2,\eta})$ for some $\eta \in \Lambda_\zeta\}$.
\end{enumerate}

So
\mn
\begin{enumerate}
\item[$\circledast_3$]  $\langle M^*_\zeta:\zeta \le \kappa\rangle$
is $\subseteq$-increasing (notice that we didn't demand
continuity) and $M^*_\kappa = M_*$.
\end{enumerate}
\mn
So it is enough to prove that $M^*_\zeta$ is
$(\lambda,\kappa)$-i.md.-limit for $\zeta < \kappa$ but not for $\zeta
= \kappa$.
\mn
\begin{enumerate}
\item[$\odot_1$]   $M^*_\kappa = M_*$ is not a
$(\lambda,\kappa)$-i.md.-limit model.
\end{enumerate}
\mn
Why?  Assume toward contradiction that there is an isomorphism $g$
from $M_\kappa$ onto $M^*_\kappa$ and let $N_i := g(M_i)$ for $i<
\kappa$, and let $h:M^*_\kappa \rightarrow \kappa$ be $h(c) 
= \min\{i < \kappa:c \in N_{i+1}\}$.  Fix $\eta \in \Lambda_\kappa$ for a while
and let $(C'_{1,\eta},C'_{2,\eta})$ be the cut of $M^*_\kappa = M_*$ with
$C'_{1,\eta} := \{c \in M_*:c <_{M_*} c_{\eta \restriction \zeta}$ for
some $\zeta < \kappa\}$.  Clearly $\langle c_{\eta \restriction
\zeta}:\zeta < \kappa\rangle$ is an increasing unbounded sequence of
members of $C'_{1,\eta}$ and $\langle f_\eta(d_\alpha):\alpha <
\lambda\rangle$ ($f_\eta$ is from $\circledast_1(c))$ is a 
decreasing sequence of members of $C'_{2,\eta}$ unbounded from below
in it.  
So $\cf(C'_{1,\eta},C'_{2,\eta}) = (\kappa,\lambda)$.
This implies that for some $i=i(\eta) < \kappa$, the set $C'_{2,\eta}
\cap N_i$ is unbounded from below in $C'_{2,\eta}$.  Hence recalling
the choice of $\bar M$ there is an
increasing continuous function $h_\eta:\kappa \rightarrow \kappa$ such
that: $\cup\{(c_{\eta \restriction h_\eta(i)},
c_{\eta \restriction j})_{M^*_\kappa}:j \in
[h_\eta(i),\kappa)\}$ is disjoint to $N_i$.  All this holds for any
  $\eta \in \Lambda_\kappa$.   Now we choose
$(\eta_\zeta,\xi_\zeta)$ by induction on $\zeta < \kappa$ such that:
\mn
\begin{enumerate}
\item[$\circledast_4$]  $(a) \quad \xi_\zeta < \kappa$ and
$\eta_\zeta \in \Lambda_{\xi_\zeta}$
\sn
\item[${{}}$]  $(b) \quad$ if $\zeta_1 < \zeta_2 < \kappa$ then
$(\eta_{\zeta_1} \restriction \xi_{\zeta_1}) \char 94 \langle 1
 \rangle \triangleleft
\eta_{\zeta_2}$ and $\xi_{\zeta_1} < \xi_{\zeta_2}$
\sn
\item[${{}}$]  $(c) \quad$ the set $\bigcup\{[c_{\eta_\zeta
\restriction \xi_\zeta},c_{\eta_\zeta \restriction \xi})_{M^*_\kappa}:\xi \in
(\xi_{\zeta +1},\kappa)\}$ is disjoint to $N_\zeta$
\sn
\item[${{}}$]  $(d) \quad$ if $\zeta$ is a successor then 
$\xi_\zeta$ is a successor
\sn
\item[${{}}$] $(e) \quad$ if $\eta_{\zeta +1 \restriction
\xi_{\zeta +1}} \triangleleft \nu \in {}^{\kappa \ge}(\lambda +1)$
then $c_\nu \in (C_{\eta_{\zeta +1 \restriction (\xi_{\zeta +1}-1) \char
94 <1>}}$,

\hskip25pt $C_{\eta_{\zeta +1} \restriction (\xi_{\zeta +1} -1)\char 94
<2>})_{M_\kappa}$ is disjoint to $N_\zeta$.
\end{enumerate}
\mn
There is no problem to carry the induction:
\medskip

\noindent
\underline{Case 1}:  $\zeta=0$.

Choose $\xi_\zeta = 0,\eta_\zeta \in \Lambda_0$.
\medskip

\noindent
\underline{Case 2}:  $\zeta = \zeta_1 +1$.

Choose $\xi_\zeta = h_{\eta_{\zeta_1}}(\xi_{\zeta_1}) + 6$.

Choose $\eta_\zeta$ such that

\[
\eta_\zeta \restriction (h_{\eta_{\zeta_1(\xi_{\zeta_1})+5)} \char 94
  \langle 1 \rangle} \trianglelefteq \eta_\zeta \in \Lambda_{\xi_\zeta}.
\]
\medskip

\noindent
\underline{Case 3}:  $\zeta$ limit.

$\xi_\zeta = \cup\{\xi_\alpha:\alpha < \zeta\}$.

Choose $\eta_\zeta \in \Lambda_{\xi_\zeta +1}$ such that $\alpha <
\zeta \Rightarrow \eta_\alpha \restriction \xi_\alpha \trianglelefteq 
\eta_\zeta$.

 Let $\eta = \bigcup\{\eta_\zeta \restriction 
\xi_\zeta:\zeta < \kappa\}$.  So $\eta
\in {}^\kappa(\lambda + 1)$ and $c_\eta \notin N_\zeta$ for every
$\zeta < \kappa$ but $\cup\{N_\zeta:\zeta < \kappa\} = M^*_\kappa =
M^*$, contradiction, so $\odot_1$ holds indeed.
\mn
\begin{enumerate}
\item[$\boxdot$]  $M^*_\zeta$ is a $(\lambda,\kappa)$-.i.md.-limit
model for $\zeta < \kappa$.
\end{enumerate}
\mn
Why?  We define $M_{\zeta,i} \subseteq M^*_\zeta$ 
for $i < \kappa$ by: $c \in M_{\zeta,i}$ iff one of the following occurs:
\mn
\begin{enumerate}
\item[$(a)$]  $c \in M_i$ but for no $\eta \in \Lambda_\zeta$ do we
have $c \in  B_\eta := \bigcup\{[c_{\eta \restriction (\zeta +i)},
c_{\eta \restriction \varepsilon})_{M_i}:\varepsilon \in (\zeta +i,\kappa)\}$
\sn
\item[$(b)$]  $c \in f_\eta(M'_i)$ for some $\eta \in \Lambda_\zeta$.
\end{enumerate}
\mn

Let
\mn
\begin{enumerate}
\item[$\bullet$]  $J_{\zeta,\eta} = 
\bigcup\{(C_{\eta \restriction \zeta},C_{\eta
\restriction \varepsilon})_{M^*_\kappa}:\varepsilon \in
(\zeta,\kappa)\}$
\sn
\item[$\bullet$]  $J_{\zeta,\eta,\varepsilon} = (C_{\eta \restriction \zeta},
C_{\eta \restriction \varepsilon})$
\sn
\item[$\bullet$]  $\langle J_{\zeta,\eta}:\eta \in \Lambda_\zeta 
\rangle \text{ are pairwise disjoint}$
\sn
\item[$\bullet$]  $J_{\zeta,\eta,\varepsilon} 
\text{ is an initial segment of } J_{\zeta,\eta}$
\sn
\item[$\bullet$]  $J_{\zeta,\eta} = 
\bigcup\{J_{\zeta,\eta,\varepsilon}:\varepsilon \in (\zeta,\kappa)\}$.
\end{enumerate}
\mn
We will make $M_{\zeta,i} \cap J_{\zeta,\eta}$ bounded in
$J_{\zeta,\eta}$ for each $i < \kappa$.

Now $\langle M_{\zeta,i}:i < \kappa\rangle$ is a
$(\lambda,\kappa)$-sequence, see Definition \ref{sl.1} hence by
\ref{sl.2}(1) the model $M^*_\zeta$ is a 
$(\lambda,\kappa)$-i.md.-limit model.
\end{PROOF}

\begin{claim}
\label{sl.21}  
If $\lambda = \lambda^{<\lambda} > \kappa = \cf(\kappa)$ 
\then \, $T_{\rd}$ has a $(\lambda,\kappa)$-limit$^+$model, 
i.e.: if $\langle M_\alpha:\alpha <
\lambda^+\rangle$ is $\subseteq$-increasing continuous sequence of
 models $\in \EC_\lambda(T_{\rd})$ and 
$M_{\alpha +1}$ is $(\lambda,\kappa)$-i.md.-limit model for every
$\alpha < \lambda^+$ \then \,: for a club of $\delta < \lambda^+$ if
$\cf(\delta) = \kappa$ then $M_\delta$ is a $(\lambda,\kappa)$-i.md.-limit.
\end{claim}

\begin{PROOF}{\ref{sl.21}}
Let $\langle M_\alpha:\alpha < \lambda^+\rangle$ be as
in the theorem and $M = \bigcup\limits_{\alpha < \lambda^+} M_\alpha$, without
loss of generality $\|M\| = \lambda^+$.  
As $\lambda = \lambda^{< \lambda}$ by \ref{sl.2}(4) we can
find a club $E$ of $\lambda^+$ such that:
\mn
\begin{enumerate}
\item[$\circledast$]   if $\alpha < \delta \in E$ and $(C_1,C_2)$
 is a cut of $M_\alpha$ of cofinality $\ne (\lambda,\lambda)$ and
some $a \in M$ realizes the cut \then \, some $a \in M_\delta$
realizes the cut.
\end{enumerate}
\mn
Let $\langle \alpha_\varepsilon:\varepsilon \le \kappa\rangle$ be an
increasing continuous sequence of ordinals from $E$ and we shall prove that
$M_{\alpha_\kappa}$ is $(\lambda,\kappa)$-i.md.-limit; this suffices
(really just $\alpha_\kappa \in E$ suffice).

Now $M_{\alpha_\kappa +1}$ is $(\lambda,\kappa)$-i.md.-limit hence
there is an $\subseteq$-increasing continuous
sequence $\langle M_{\alpha_\kappa
+1,i}:i < \kappa \rangle$ witnessing $M_{\alpha_\kappa +1}$ is
$(\lambda,\kappa)$-i.md.-limit model, i.e. its union is
$M_{\alpha_\kappa +1}$ and it is a $(\lambda,\kappa)$-sequence.  
Now $M_{\alpha_\kappa +1,i} \cap
M_{\alpha_\kappa} = \bigcup\{M_{\alpha_{\kappa +1,i}} \cap
M_{\alpha_\zeta}:\zeta < \kappa\}$ but $\kappa < \lambda = \cf(\lambda)$
 hence \wilog \, $M_{\alpha_\kappa,0}
\cap M_{\alpha_0}$ has cardinality $\lambda$ hence $N_i :=
M_{\alpha_\kappa +1,i} \cap M_{\alpha_i} \in \EC_\lambda(T_{\rd})$.

Clearly
\mn
\begin{enumerate}
\item[$(*)_1$]   $\langle N_i:i < \kappa\rangle$ is a
$\subseteq$-increasing continuous sequence of members of $\EC_\lambda(T_{\rd})$ 
with union $M_{\alpha_\kappa}$.
\end{enumerate}
\mn
So it is enough to show that $\langle N_i:i < \kappa\rangle$ is a
$(\lambda,\kappa)$-sequence by \ref{sl.2}(1).  By
$(*)_1$, clause (a) from Definition \ref{sl.1} holds.
\mn
\begin{enumerate}
\item[$\circledast_2$]  $\langle N_i:i < \kappa\rangle$ satisfies
clause (b) of \ref{sl.1}.
\end{enumerate}
\mn
[Why?  Let $i < \kappa$ and $(C_1,C_2)$ be a cut of $N_i$.  First, we
assume $(C_1,C_2)$ is of cofinality $\ne (\lambda,\lambda)$.  
As $C_1,C_2 \subseteq N_i \subseteq M_{\alpha_\kappa +1,i}$ by
the properties of $\langle M_{\alpha_\kappa +1,i}:i < \kappa\rangle$
there is $a \in M_{\alpha_\kappa +1}$ such that $C_1 < a < C_2$.
  
If for some $b \in M_{\alpha_i},C_1 < b < C_2$ then \wilog \, $a \in
M_{\alpha_\kappa}$ and we are done.  If not, $a$ induces on
$M_{\alpha_i}$ a cut $(C'_1,C'_2),C_1 \subseteq C'_1,C_2 \subseteq
C'_2$, cf$(C'_1,C'_2) = \cf(C_1,C_2) \ne (\lambda,\lambda)$.
As $\alpha_i < \alpha_\kappa \in E$, by $\circledast$
there is $a \in M_{\alpha_i+1} \subseteq 
M_{\alpha_\kappa}$ such that $C_1 < a < C_2$.  So
clause (b) of Definition \ref{sl.1} really holds.

Second, we assume that $(C_1,C_1)$ is of cofinality
$(\lambda,\lambda),\kappa < \lambda = \cf(\lambda)$ so \wilog \,
clause $(b)(\beta)$ of \ref{sl.1} holds so some $a \in M_{\alpha_\kappa}$
realizes $(C_1,C_2)$ so for some $j \in (i,\kappa),a \in M_{\alpha_j}$
hence the cut $(\{b \in M_{\alpha_j}:b < a\},\{c \in M_{\alpha_j}:a
\le c\})$ of $M_{\alpha_j}$ has cofinality $\ne (\lambda,\lambda)$ so
is realized by infinitely many $a' \in M_{\alpha_j+1} \subseteq
M_{\alpha_\kappa}$, hence also $(C_1,C_2)$ is, so clause $(b)(\beta)$
of \ref{sl.1} holds.  Together $(C_1,C_1)$ satisfies $(\alpha)$ of
$(\beta)$ of \ref{sl.1} as promised.]
\mn
\begin{enumerate}
\item[$\circledast_3$]   if $a <_{M_{\alpha_\kappa}} b$ then
$M_{\alpha_\kappa} \restriction (a,b)$ is universal (for
($\EC_{T_{\oor}}(\lambda),\subseteq))$.
\end{enumerate}
\mn
[Why?  As $\langle \alpha_\varepsilon:\varepsilon \le \kappa\rangle$
is increasing continuous and $\langle M_\alpha:\alpha <
\lambda^+\rangle$ is increasing continuous, clearly
for some $i < \kappa$ we have $a,b \in
M_{\alpha_i}$ hence $M_{\alpha_i+1} \restriction (a,b)$ is $\lambda$- universal
but $M_{\alpha_i+1} \subseteq M_{\alpha_k}$ so $M_{\alpha_\kappa}
\restriction (a,b)$ is universal so we are done.]
\end{PROOF}

\bibliographystyle{plain}
\bibliography{lista,listb,listx,listf,liste,listy}

\end{document}